\documentclass{article}
\usepackage{amssymb}
\usepackage{graphicx}
\usepackage{amsmath}

\newtheorem{theorem}{Theorem}
\newtheorem{acknowledgement}[theorem]{Acknowledgement}

\newtheorem{corollary}[theorem]{Corollary}

\newtheorem{definition}[theorem]{Definition}
\newtheorem{example}[theorem]{Example}

\newtheorem{lemma}[theorem]{Lemma}

\newtheorem{proposition}[theorem]{Proposition}
\newtheorem{remark}[theorem]{Remark}

\newenvironment{proof}[1][Proof]{\textbf{#1.} }{\ \rule{0.5em}{0.5em}}
\input{tcilatex}

\begin{document}

\title{Approximations of the Brownian Rough Path with Applications to
Stochastic Analysis.\\
Des Approximations du Rough Path Brownien et Applications \`a l'Analyse
Stochastique.}
\author{Peter Friz \thanks{%
Courant Institute, 251 Mercer St. New York, NY 10012, USA,
Peter.Friz@cims.nyu.edu} and Nicolas Victoir\thanks{%
Magdalen College, Oxford OX1 4AU, UK; victoir@maths.ox.ac.uk}}
\maketitle

\begin{abstract}
A geometric $p$-rough path can be seen to be a genuine path of finite $p$%
-variation with values in a Lie group equipped with a natural distance. The
group and its distance lift $(\mathbb{R}^{d},+,0)$ and its Euclidean
distance.

This approach allows us to easily get a precise modulus of continuity for
the Enhanced Brownian Motion (the Brownian Motion and its Levy Area).

As a first application, extending an idea due to Millet \&\ Sanz-Sole, we
characterize the support of the Enhanced Brownian Motion (without relying on
correlation inequalities). Secondly, we prove Schilder's theorem for this
Enhanced Brownian Motion. As all results apply in H\"{o}lder (and stronger)
\ topologies, this extends recent work by Ledoux, Qian, Zhang \cite{LQZ}.
Lyons' fine estimates in terms of control functions \cite{Ly} allow us to
show that the It\^{o} map is still continuous in the topologies we
introduced. This provides new and simplified proofs of the Stroock-Varadhan
support theorem and the Freidlin-Wentzell theory. It also provides a short
proof of modulus of continuity for diffusion processes along old results by
Baldi.
\end{abstract}

\begin{abstract}
Un $p$-rough path est un chemin de $p$-variation finie \`{a} valeurs dans un
groupe de Lie muni d'une distance sous-riemannienne. Le groupe et sa
distance g\'{e}neralisent $(\mathbb{R}^{d},+,0)$ et la distance Euclidienne.%
\newline
Cette approche nous permet d'obtenir un modulus de continuit\'{e} tres pr%
\'{e}cis pour le rough path brownien (le mouvement brownien et son aire de
Levy). Pour ce dernier, nous prouvons un th\'{e}or\`{e}me du support
(adaptant une id\'{e}e de Millet et Sans-Sole) et un th\'{e}or\`{e}me de
Schilder. Comme tous les r\'{e}sultats sont prouv\'{e}s en utilisant des
topologies de type H\"{o}lder ou plus fines, cela g\'{e}neralise le papier
de Ledoux, Qian, Zhang \cite{LQZ}. Les r\'{e}sultats de T.Lyons \cite{Ly}
permettent de prouver rapidement que la fonction d'It\^{o} est continue pour
les topologies que nous avons introduites. Cela nous donne des nouvelles
preuves du th\'eor\`eme du support de Stroock-Varadhan et de la th\'eorie de
Freidlin-Wentzell. Nous obtenons au passage une preuve simple du modulus de
continuit\'e pour les processus de diffusions, obtenu pr\'{e}c\'{e}demment
par Baldi.
\end{abstract}

\section{Introduction}

\bigskip Starting with \cite{L94}, Terry Lyons developed a general theory of
integration and differential equations of the form 
\begin{equation}
dy_{t}=f(y_{t})dx_{t}.  \label{ODE}
\end{equation}
To include the important example of stochastic differential equations, $x$
is allowed to be ``rough'' in some sense. Standard H\"{o}lder regularity of
Brownian motion, for instance, implies finite $p$-variation only for $p>2$.
Another issue was to explain (deterministically) the difference between
stochastic differential equations based on Stratonovich versus It\^{o}
integration. Last but not least, motivated from examples like Fractional
Brownian motion, driving signals much rougher than Brownian motion should be
included. \newline
All this has been accomplished in a beautiful way and the reader can
nowadays find the general theory exposed in ~\cite{Lej,Ly,LQ}.

Loosely speaking, for general $p\geq 1$, one needs to \textquotedblleft
enhance\textquotedblright\ the driving signal $x$, with values in some
Banach space $V$, to $X\in V\oplus V^{\otimes 2}...\oplus V^{\otimes \lbrack
p]}$ such that the resulting object $X$ satisfies certain algebraic 
\footnote{%
For \textit{algebraic} convenience $X$ is often enhanced to $\mathbb{R}%
\oplus V\oplus V^{\otimes 2}...\oplus V^{\otimes \lbrack p]}$ with scalar
component constant $1$.} and analytic conditions. For $x$ of finite
variation, this enhancement will simply consist of all the iterated
integrals of $x$, 
\begin{equation*}
X_{s,t}^{k}:=\int_{s<u_{1}<...<u_{k}<t}dx_{u_{1}}\otimes ...\otimes
dx_{u_{k}},\ \ \ k=1,...,[p].
\end{equation*}
These are the \textit{Smooth Rough Paths.} Consider a time horizon of $[0,1]$
(valid for the rest of the paper) and introduce the $p$-variation metric,
defined as 
\begin{equation*}
d(X,Y)=\max_{k=1,...,[p]}\left(
\sup_{D}\sum_{l}|X_{t_{l-1},t_{l}}^{k}-Y_{t_{l-1},t_{l}}^{k}|^{p/k}\right)
^{k/p},
\end{equation*}
where $\sup_{D}$ runs over all finite divisions of $[0,1]$. Here $|.|$
denotes (compatible) tensor norms in $V^{\otimes k}$. Closure of Smooth
Rough Paths with respect to this metric yields the class of \textit{%
Geometric Rough Paths in the sense of }\cite{LQ}, denoted by $G\Omega
_{p}(V) $. The solution map, also called \textit{It\^{o} map}, to (\ref{ODE}%
) is then a continuous map from $G\Omega _{p}(V)\rightarrow G\Omega _{p}(W)$%
, provided $f:W\rightarrow L(V,W)$ satisfies mild regularity conditions.
This is Lyons' celebrated \textit{Universal Limit Theorem}. In particular,
smooth approximations $X(n)$ which converge in $p$-variation to $X\in
G\Omega _{p}(V)$ will cause the corresponding solutions $Y(n)$ to converge
to $Y$ in $p$-variation. Hence, one deals with some kind of generalized
Stratonovich theory.

\textit{However}, the so important case of $p\in (2,3)$, on which this paper
will focus, allows for more. For the sake of concreteness, we will set $V=%
\mathbb{R}^{d}$ from here on. Following \cite[p. 149]{LQ} and also~\cite{Ly}
the driving signal only needs to be a \textit{Multiplicative functional of
finite $p$-variation}. By definition, this is a continuous map 
\begin{equation*}
(s,t)\rightarrow (X_{s,t}^{1},X_{s,t}^{2})\in \mathbb{R}^{d}\oplus \left( 
\mathbb{R}^{d}\right) ^{\otimes 2}=:T^{2},
\end{equation*}
where $0\leq s\leq t\leq 1$, satisfying the algebraic \textit{Chen condition 
} 
\begin{equation}
X_{s,u}=X_{s,t}\otimes X_{t,u}\Leftrightarrow \text{ }%
X_{s,u}^{1}=X_{s,t}^{1}+X_{t,u}^{1},\text{ }%
X_{s,u}^{2}=X_{s,t}^{2}+X_{t,u}^{2}+X_{s,t}^{1}\otimes X_{t,u}^{1},
\label{chen}
\end{equation}
whenever $s\leq t\leq u$, and the analytic condition $d_{p\text{-var}%
}(X,0)<\infty $ \ i.e. 
\begin{equation}
\sup_{(0\leq t_{0}<\ldots <t_{n}\leq
1)}\sum_{l}|X_{t_{l-1},t_{l}}^{k}|^{p/k}<\infty ,\ \ \ \ k=1,2.
\label{finite_p_var}
\end{equation}
(Often $k=1,2$ are referred to as \textit{first} resp. \textit{second level}%
). The class of such rough paths is denoted $\Omega _{p}(\mathbb{R}^{d})$.
Condition (\ref{chen}) is known as \textit{Chen relation} and expresses
simple additive properties whenever $X^{2}$ is obtained as \textit{some}
iterated integral including the cases of Stratonovich resp. It\^{o} Enhanced
Brownian Motion. Whenever a first order calculus underlies this integration
(which is the case for Stratonovich integration), one has 
\begin{equation}
Symm(X^{2})=\frac{1}{2}X^{1}\otimes X^{1}  \label{geometric}
\end{equation}

and such paths are called \textit{Geometric Rough Paths of finite }$p$%
\textit{-variation\ }or simply \textit{Geometric }$p$\textit{-Rough Paths, }%
write\textit{\ }$X\in G\Omega (\mathbb{R}^{d})^{p}.$ Clearly,\textit{\ } 
\begin{equation*}
\{\text{Smooth Rough Paths}\}\subset G\Omega _{p}(\mathbb{R}^{d})\subset
G\Omega (\mathbb{R}^{d})^{p}\subset \Omega _{p}(\mathbb{R}^{d}).
\end{equation*}

\bigskip One can indeed choose in which space to work with and the Lyons
theory will provide meaning, existence and uniqueness to the purely
deterministic \textit{rough differential equation} 
\begin{equation*}
dY=f(y_{0}+Y_{0t}^{1})dX
\end{equation*}

where $f=(V_{1},...,V_{d})$ are, in general non-commuting, vector fields
with mild regularity conditions. As before, the It\^{o} map $X\mapsto Y$ is,
continuous under $p$-variation topology. This rough differential equation
indeed generalizes ordinary and stochastic (Stratonovich and It\^{o})
differential equations. For instance, it is known (and also follows from the
results in this paper) that a.s. the Stratonovich Enhanced Brownian Motion
(EBM) $\mathbf{B}\in G\Omega _{p}(\mathbb{R}^{d})$ for all $p\in (2,3)$.
Choosing $X=\mathbf{B}$ the projection of the rough path $Y$ to its first
level will solve the associated Stratonovich stochastic differential
equation. That is 
\begin{equation*}
y_{0}+Y_{0t}^{1}\text{ solves }dy=\sum_{i}V_{i}(y)\circ d\beta ^{i}.
\end{equation*}

\bigskip

For the rest of the paper, $p$ denotes a fixed real in $(2,3)$. The
contributions of this paper may be summarized as follows:

\bigskip

(a) We look at geometric $p$-rough paths from a new angle. Observe that $%
G:=\{X\in T^{2}:(\ref{geometric})$ holds$\}$ is the free nilpotent Lie group
of step 2 \cite{Ly,Re}, a simply connected Lie group which lifts $(\mathbb{R}%
^{d},+,0)$.\ Chen's condition is equivalent to the fact that $\mathbf{x}%
_{t}=X_{0,t}$ is a $G$-valued path such that $X_{s,t}=\mathbf{x}%
_{s}^{-1}\otimes \mathbf{x}_{t}=\mathbf{x}_{s,t}.$ We put a homogenous,
sub-additive norm on $\left( G,\otimes \right) $. Geometric $p$-rough paths
are then easily seen to be $G$-valued paths of finite $p$-variation.
Standard proofs for Kolmogorov's criterion or the Garsia, Rumsey, Rodemich
inequality adapt with no changes from $(\mathbb{R}^{d},+)$-valued to $%
(G,\otimes )$-valued processes. With this observation, regularity results
for the EBM $\mathbf{B},$ as H\"{o}lder continuity and L\'{e}vy modulus of
continuity, follow after simple moment estimates. Sometimes, it will be
convenient to work in the associated\ Lie algebra of the group $G$. For
instance, the EBM\ viewed through this chart is nothing else than the well
studied Gaveau diffusion \cite{Ga}. Its importance in the context of limit
theorems was already highlighted in Malliavin's book, \cite{Ma}.

\bigskip

(b)\ We introduce a number of different topologies on $G\Omega _{p}(\mathbb{R%
}^{d}),$ effectively reducing this space to geometric rough paths for which
the associated norms are finite. For instance, we are able to deal with H%
\"{o}lder and ``modulus type'' norms. Exploiting fine estimates in Lyons'
Limit Theorem we have continuity of the It\^{o} map in all these topologies.
In former applications of rough path theory to stochastic analysis result
were always obtained in $p$-variation topology, leaving open a gap between
(usually well known) results in H\"{o}lder and stronger topologies.

\bigskip

(c) Lyons' Universal Limit Theorem implies a L\'{e}vy modulus of continuity
for diffusions, along the results by Baldi \cite{Ba1}.

\bigskip

(d) We establish convergence of several different approximations to the EBM.
The EBM is usually defined as the limit of a sequence of smooth rough paths,
which is shown to be Cauchy \cite{LQ}. Here, we define directly the EBM, and
its regularity allows us to prove convergence of some sequences of smooth
rough paths to the EBM. The first idea, common to works in \cite{Le,Ma,HL},
is that approximations are obtained by conditioning with respect to dyadic
filtrations. After establishing uniform regularity of approximations, easily
obtained by Doob's inequality, we can use a compactness argument to show
convergence in interesting topologies of some sequences of smooth rough
paths to the EBM.

\bigskip

(e) By combining ideas due to Millet and Sanz-Sole with our results above we
give a short and original proof of the support theorem for EBM\ in the
strong topologies mentioned above. By means of \ continuity of the It\^{o}
map, from (b), this immediately implies the support theorem for diffusions.
In the context of rough paths this improves work by Ledoux, Qian, Zhang
(support theorem in $p$-variation topology,\cite{LQZ}) and by the first
cited author (H\"{o}lder topology,\cite{F}).

As for the history of the support theorem, it was originally obtained by
Stroock, Varadhan \cite{SV2} in sup topology, then by Ben Arous, Gradinaru,
Ledoux \cite{BA2} and Millet, Sanz-Sole \cite{Millet} in H\"{o}lder norm of
exponent less than $1/2.$ Extension to modulus space has been obtain in \cite%
{GNS} and to Orlicz-Besov space in \cite{Me}. We limit ourselves to
``modulus norm'', in the spirit of \cite{GNS}, and we will not recover fully
the results in \cite{GNS,Me}. On the other hand, we have a description of
the support of the Stratonovich enhanced diffusion the goes beyond the last
quoted results.

\bigskip

(f) Schilder's theorem for EBM is obtained. As before, continuity of the It%
\^{o} map will give the Freidlin-Wentzell Large Deviation result in the
topologies mentioned in (b). Strassen's law is obtained as corollary. Again,
we improve \cite{LQZ} and recover well known large deviations results in H%
\"{o}lder and modulus norm \cite{BL}. Once again though, we do not deal with
Orlicz-Besov metrics.

\bigskip

Constants in this paper may varies from line to line.

\begin{acknowledgement}
\bigskip The authors would like to thank G. Ben Arous, T. Lyons and S.
Varadhan for related discussions.
\end{acknowledgement}

\section{Rough paths}

\subsection{Free Nilpotent Lie Group of step 2}

We fix the dimension $d$ ($d\geq 2$ to avoid trivialities) and we denote by $%
\mathcal{L}\left( \mathbb{R}^{d}\right) =\mathbb{R}^{d}\oplus so(d)$, where $%
so(d)$ $\mathbb{\simeq R}^{d(d-1)/2}$ is the space of real antisymmetric $%
d\times d$ matrices. With the bracket 
\begin{eqnarray*}
\lbrack ,]:\mathcal{L}\left( \mathbb{R}^{d}\right) \times \mathcal{L}\left( 
\mathbb{R}^{d}\right) &\rightarrow &\mathcal{L}\left( \mathbb{R}^{d}\right)
\\
\left( \left( a^{1},a^{2}\right) ,\left( b^{1},b^{2}\right) \right)
&\rightarrow &[a,b]=(0,a^{1}\otimes b^{1}-b^{1}\otimes a^{1}),
\end{eqnarray*}
$\mathcal{L}\left( \mathbb{R}^{d}\right) $ becomes a (step 2 nilpotent) Lie
algebra. The group multiplication in the associated simply connected Lie
Group $G(\mathbb{R}^{d})=\exp (\mathcal{L}\left( \mathbb{R}^{d}\right) ),$%
\cite{Re,Var,Wa}, is given by the Baker-Campbell-Hausdorff formula\footnote{%
Due to step 2 nilpotency, only the first bracket appears.}.

\begin{eqnarray*}
\otimes :G\left( \mathbb{R}^{d}\right) \times G\left( \mathbb{R}^{d}\right)
&\rightarrow &G\left( \mathbb{R}^{d}\right) \\
\exp (a)\otimes \exp (b) &\rightarrow &\exp \left( a+b+\frac{1}{2}%
[a,b]\right) .
\end{eqnarray*}
Its neutral element is $\exp (0)$, and the inverse of $\exp (a)$ is $\exp
(-a)$. We can identify $G\left( \mathbb{R}^{d}\right) $ with the nonlinear
submanifold of $\mathbb{R}^{d}\oplus \mathbb{R}^{d\times d}$ given by 
\begin{equation*}
\{g=(g^{1},g^{2})\in \mathbb{R}^{d}\oplus (\mathbb{R}^{d})^{\otimes 2}:\text{%
symmetric part of }g^{2}\text{ equals }\frac{1}{2}g^{1}\otimes g^{1}\}
\end{equation*}
with usual (truncated)\ tensor multiplication, that is, 
\begin{equation*}
g\otimes h=(g^{1}+h^{1},g^{2}+g^{1}\otimes h^{1}+h^{2}).
\end{equation*}

$\left( G\left( \mathbb{R}^{d}\right) ,\otimes ,\exp (0)\right) $ is the
free nilpotent group of step $2$ over $\mathbb{R}^{d}$, \cite{Ly,Re}. Note
that $G\left( \mathbb{R}^{d}\right) $ is invariant under the dilation
operator $\delta _{t}$, for $t\in \mathbb{R}$, $\delta _{t}$ being defined
by 
\begin{eqnarray*}
G\left( \mathbb{R}^{d}\right) &\rightarrow &G\left( \mathbb{R}^{d}\right) \\
\exp (a^{1},a^{2}) &\rightarrow &\exp (ta^{1},t^{2}a^{2}).
\end{eqnarray*}%
We define on the group 
\begin{equation}
||g||\text{ }=\inf_{\substack{ x^{1},\ldots ,x^{n}\in \mathbb{R}^{d}  \\ %
\bigotimes_{i=1}^{n}\exp (x^{i})=g}}\sum_{i=1}^{n}\left\vert
x_{i}\right\vert _{\mathbb{R}^{d}},  \label{def_of_norm_on_G}
\end{equation}%
where $\left\vert .\right\vert _{\mathbb{R}^{d}}$ is the Euclidean norm on $%
\mathbb{R}^{d}$. $\left\Vert .\right\Vert $ is a sub-additive, symmetric
homogeneous norm \cite{FS} on $G\left( \mathbb{R}^{d}\right) $ , that is%
\newline
\begin{tabular}{ll}
(i) & $\left\Vert g\right\Vert $ if and only if $g=\exp (0)$, \\ 
(ii) & for all $g\in G\left( \mathbb{R}^{d}\right) $ and $t\in \mathbb{R}$, $%
\left\Vert \delta _{t}g\right\Vert =\left\vert t\right\vert \left\Vert
g\right\Vert ,$ \\ 
(iii) & for all $g,h\in G\left( \mathbb{R}^{d}\right) $, $\left\Vert
g\otimes h\right\Vert \leq \left\Vert g\right\Vert +\left\Vert h\right\Vert
, $ \\ 
(iv) & for all $g$, $\left\Vert g\right\Vert =\left\Vert g^{-1}\right\Vert .$%
\end{tabular}
\newline
From this sub-additive, symmetric homogeneous norm, we construct a left
invariant distance on $G\left( \mathbb{R}^{d}\right) $ (which is a
Carnot-Caratheodory distance \cite{Gr,Mo}) by defining 
\begin{equation*}
d(g,h)=\left\Vert h^{-1}\otimes g\right\Vert .
\end{equation*}%
If $|.|_{\mathbb{R}^{d}\otimes \mathbb{R}^{d}}$ denotes a norm on $\mathbb{R}%
^{d}\otimes \mathbb{R}^{d}$, then

\begin{equation*}
|||\exp (a^{1},a^{2})|||\text{ }=|a^{1}|_{\mathbb{R}^{d}}+\sqrt{|a^{2}|_{%
\mathbb{R}^{d}\otimes \mathbb{R}^{d}}}
\end{equation*}
defines another homogeneous norm on $G\left( \mathbb{R}^{d}\right) $ (that
is a norm satisfying (i) and (ii)), and as all homogeneous norms are
equivalent \cite{HS}, one can find some positive constants $c_{1},c_{2}$
such that for all $g\in G\left( \mathbb{R}^{d}\right) $ 
\begin{equation}
c_{1}|||g|||\leq \left\| g\right\| \leq c_{2}|||g|||.  \label{boundnorm}
\end{equation}
This implies the following:

\begin{corollary}
For some constant $C$, 
\begin{equation}
\left\Vert h^{-1}\otimes g\otimes h\right\Vert \leq C\left( \left\Vert
g\right\Vert +\sqrt{\left\Vert h\right\Vert \left\Vert g\right\Vert }\right)
,  \label{useful_ineq}
\end{equation}%
and for any $k\geq 2,$ 
\begin{equation}
d\left( \bigotimes_{i=1}^{k}g_{i},\bigotimes_{i=1}^{k}h_{i}\right) \leq
C\sum_{i=1}^{k}\left( d(g_{i},h_{i})+\sqrt{d(g_{i},h_{i})\left\Vert
\bigotimes_{j=i+1}^{k}h_{j}\right\Vert }\right) .  \label{useful_ineq2}
\end{equation}
\end{corollary}

\begin{proof}
If $g=e^{b},h=e^{a}$, $h^{-1}\otimes g\otimes h=e^{b}\otimes e^{[a,b]}$,
hence 
\begin{eqnarray*}
\left\vert \left\vert \left\vert h^{-1}\otimes g\otimes h\right\vert
\right\vert \right\vert &\leq &\left\vert \left\vert \left\vert g\right\vert
\right\vert \right\vert +\left\vert \left\vert \left\vert
e^{[a,b]}\right\vert \right\vert \right\vert \\
&\leq &\left\vert \left\vert \left\vert g\right\vert \right\vert \right\vert
+\sqrt{\left\vert \left\vert \left\vert g\right\vert \right\vert \right\vert
.\left\vert \left\vert \left\vert h\right\vert \right\vert \right\vert }.
\end{eqnarray*}%
The inequality is then proved using inequality (\ref{boundnorm}). The second
inequality is a consequence of the first one. We show it for $k=2$, the
general case follows in exactly the same way, by induction.%
\begin{eqnarray*}
d\left( g_{1}\otimes g_{2},h_{1}\otimes h_{2}\right) &=&\left\Vert
h_{2}^{-1}\otimes h_{1}^{-1}\otimes g_{1}\otimes g_{2}\right\Vert \\
&=&\left\Vert h_{2}^{-1}\otimes h_{1}^{-1}\otimes g_{1}\otimes h_{2}\otimes
h_{2}^{-1}\otimes g_{2}\right\Vert \\
&\leq &\left\Vert h_{2}^{-1}\otimes h_{1}^{-1}\otimes g_{1}\otimes
h_{2}\right\Vert +\left\Vert h_{2}^{-1}\otimes g_{2}\right\Vert \\
&\leq &C\left( d\left( g_{1},h_{1}\right) +\sqrt{d\left( g_{1},h_{1}\right)
\left\Vert h_{2}\right\Vert }\right) +d\left( g_{2},h_{2}\right) .
\end{eqnarray*}
\end{proof}

$\left( G\left( \mathbb{R}^{d}\right) ,\otimes ,\exp (0)\right) $ equipped
with a homogeneous norm is a simple generalization of $\left( \mathbb{R}%
^{d},+,0\right) $ equipped with a norm.

We let $C_{0}\left( [0,1],G\left( \mathbb{R}^{d}\right) \right) $ to be the
space of continuous function from $[0,1]$ to $G\left( \mathbb{R}^{d}\right) $
such that their value at time $0$ is $\exp (0)$. With a slight abuse, we
will call such elements $G\left( \mathbb{R}^{d}\right) $-valued paths. If $%
\mathbf{x}\in C_{0}\left( [0,1],G\left( \mathbb{R}^{d}\right) \right) $ and $%
s<t$, we will denote by $\mathbf{x}_{s,t}$ the element $\mathbf{x}%
_{s}^{-1}\otimes \mathbf{x}_{t}$.

\begin{remark}
Let $|.|_{\mathbb{R}^{d}\otimes \mathbb{R}^{d}}$ be a compatible tensor norm%
\footnote{%
That is $|a\otimes b|_{\mathbb{R}^{d}\otimes \mathbb{R}^{d}}$ $\leq |a|_{%
\mathbb{R}^{d}}|b|_{\mathbb{R}^{d}}.$}$,$ an explicit norm satisfying
(i)-(iv) is given by $\max \left\{ |a^{1}|_{\mathbb{R}^{d}},\sqrt{|a^{2}+%
\frac{1}{2}a^{1}\otimes a^{1}|_{\mathbb{R}^{d}\otimes \mathbb{R}^{d}}}%
\right\} $ where $(a^{1},a^{2})\in \mathcal{L}\left( \mathbb{R}^{d}\right) .$
\end{remark}

\subsection{$p$-Variation}

Let $\left( G,\otimes ,e\right) $ be a group equipped with a homogeneous
norm $\left\| .\right\| .$ Here, we think of $G$ being either $\left( 
\mathbb{R}^{d},+,0\right) $ or $\left( G\left( \mathbb{R}^{d}\right)
,\otimes ,\exp (0)\right) .$ A path $x:[0,1]\rightarrow G$ is said to have
finite $p$-variation if 
\begin{equation*}
\sup_{(0\leq t_{0}<\ldots <t_{n}\leq 1)}\sum_{i}\left\|
x_{t_{i},t_{i+1}}\right\| ^{p}<\infty ,
\end{equation*}
Note that a path $x$ is continuous and of finite $p$-variation if and only
if (see \cite{LQ}) 
\begin{equation*}
\left\| x_{s,t}\right\| ^{p}\leq \omega (s,t)\text{ \ \ \ \ \ for all }s\leq
t
\end{equation*}
for some control function $\omega .$ By definition, this means 
\begin{equation}
\begin{tabular}{ll}
(i) & $\omega :\left\{ (s,t),0\leq s\leq t\leq 1\right\} \rightarrow \mathbb{%
R}^{+}$ is continuous near the diagonal. \\ 
(ii) & $\omega $ is super-additive, i.e. $\forall $ $s<t<u$, $\omega
(s,t)+\omega (t,u)\leq \omega (t,u)$. \\ 
(iii) & $\omega (t,t)=0$ for all $t\in \lbrack 0,1]$%
\end{tabular}
\label{control}
\end{equation}
\newline
We will say in such case that $x$ has finite $p$-variation controlled by $%
\omega $. We will construct control functions in the following way,

\begin{proposition}
Consider a continuous map $f:\mathbb{R}^{+}\longrightarrow \mathbb{R}^{+},$
increasing, convex, $\ f(0)=0$. Then $(s,t)\longmapsto f(t-s)$ is a control
function.

\begin{example}
$f(t-s)=c(t-s)$ for some constant $c>0.$ This is equivalent to $1/p$-H\"{o}%
lder continuity for the controlled path.
\end{example}
\end{proposition}

\subsection{Definition of a Geometric $p$-Rough Path}

We will denote by $\pi _{1}$ resp. $\pi _{2}$ the natural injection from $%
G\left( \mathbb{R}^{d}\right) $ onto $\mathbb{R}^{d}$ resp. onto $so(d)$. If 
$\mathbf{x}$ is a $G\left( \mathbb{R}^{d}\right) $-valued path of finite $p$%
-variation, then clearly, $\pi _{1}(\mathbf{x):}t\rightarrow \pi _{1}(%
\mathbf{x}_{t})$ is a $\mathbb{R}^{d}$-valued path of finite $p$-variation.
We will say that $\mathbf{x}$ lies above $\pi _{1}(\mathbf{x)}$. Conversely,
assume $x$ is a $\mathbb{R}_{d}$-valued path of finite $p$-variation. If $%
p<2 $, then 
\begin{equation}
S(x\mathbf{):}t\rightarrow \exp \left( x_{t}+\frac{1}{2}\int_{0}^{t}\left(
x_{u}\otimes dx_{u}-dx_{u}\otimes x_{u}\right) \right)  \label{signature}
\end{equation}%
is the unique $G\left( \mathbb{R}^{d}\right) $-valued path of finite $p$%
-variation lying above $x$ ( $\pi _{2}\left( S(x\mathbf{)}\right) _{t}$ is
the Levy area of $x$ between time $0$ and $t$; such integrals are well
defined Young integrals \cite{Yo}). The path $t\rightarrow S(x)_{t}$ is
called the signature of the path $x$. \newline
If $p\in (2,3),$ then there exists a $G\left( \mathbb{R}^{d}\right) $-valued
path $\mathbf{x}$ lying above $x$ \cite{Vi}, but uniqueness is not true
anymore \cite{Ly,Vi}.

\begin{definition}
For $p\in \lbrack 2,3),$ a geometric $p$-rough path is a $G\left( \mathbb{R}%
^{d}\right) $-valued path of finite $p$-variation. Such class is denoted by $%
G\Omega \left( \mathbb{R}^{d}\right) ^{p}$.
\end{definition}

Note that if $\mathbf{x}$ is a $G\left( \mathbb{R}^{d}\right) $ valued path,
then $(s,t)\rightarrow \mathbf{x}_{s}^{-1}\otimes \mathbf{x}_{t}$ is a
geometric multiplicative functional, in the sense of \cite{Ly}.
Reciprocally, if $\mathbf{y}_{s,t}$ is a multiplicative functional, then $%
\mathbf{x}_{t}=\mathbf{y}_{0,t}$ is a $G\left( \mathbb{R}^{d}\right) $ path
starting at $\exp (0)$ and $\mathbf{y}_{s,t}=\mathbf{x}_{s}^{-1}\otimes 
\mathbf{x}_{t}.$ Geometric multiplicative functionals and $G\left( \mathbb{R}%
^{d}\right) $-valued path are the same objects. The main theorem in rough
path theory is the continuity of the It\^{o} map.

\subsection{The It\^{o} Map}

\begin{theorem}
\label{ItoMap}Let $\mathbf{x}\in G\Omega (\mathbb{R}^{d})^{p}$, $\epsilon >0$
and $\mathcal{V}$ be a linear map from $\mathbb{R}^{d}$ into the $%
Lip[p+\varepsilon ,\mathbb{R}^{n}]$ vector fields\footnote{%
A function $f$ which is $(k+\epsilon )$-Lipschitz on $\mathbb{R}^{d},$ for $%
k\in \mathbb{N}$ is a $k$ times differentiable function whose $k$-th
derivative is $\epsilon $-H\"{o}lder, using the classical definition of H%
\"{o}lder functions with parameter in $[0,1).$ See \cite{Ly,Ste}}. There
exists a solution $\mathbf{y}\in G\Omega (\mathbb{R}^{n})^{p}$ to the rough
differential equation 
\begin{equation}
\begin{tabular}{ll}
$d\mathbf{y}_{t}\mathbf{=}\mathcal{V}(\pi _{1}\left( \mathbf{y}_{t}\right) )d%
\mathbf{x}_{t},$ & $y_{0}^{1}=a$%
\end{tabular}
\ \ ,  \label{DE}
\end{equation}
i.e. there exists an extension of $\mathbf{x}$ to $\mathbf{z}\in G\Omega (%
\mathbb{R}^{d}\oplus \mathbb{R}^{n})^{p}$ such that $\mathbf{z}$ projects
onto 
\begin{equation*}
\mathbf{z}_{t}=(\mathbf{x}_{t},\mathbf{y}_{t}),
\end{equation*}
and $\mathbf{z}$ satisfies $\mathbf{z}=\int h\left( \pi _{1}\left( \mathbf{z}%
\right) \right) \delta \mathbf{z,}$ with 
\begin{eqnarray*}
h &:&\mathbb{R}^{d}\oplus \mathbb{R}^{n}\rightarrow Hom(\mathbb{R}^{d}\oplus 
\mathbb{R}^{n},\mathbb{R}^{d}\oplus \mathbb{R}^{n}) \\
(x,y) &\rightarrow &\left( \left( dX,dY\right) \rightarrow \left(
dX,f(y)dY\right) \right) .
\end{eqnarray*}
If the $p$-variation of $\mathbf{x}$ is controlled by $\omega $, then the $p$%
-variation of $\mathbf{z}$ (and hence $\mathbf{y}$) is controlled by $%
C\omega $, where $C$ depends on $p,\varepsilon ,$the H\"{o}lder norm of $%
\mathcal{V}$ and the supremum of $\omega $ on the consider interval.
Moreover, for all $s<t$ such that $\omega (s,t)<1$, 
\begin{equation*}
\left\| \mathbf{z}_{s,t}\right\| ^{p}\leq K_{p,\varepsilon ,f}\omega (s,t),
\end{equation*}
where $K_{p,\varepsilon ,\mathcal{V}}$ is a constant which only depends on $%
p $,$\varepsilon $, and the H\"{o}lder norm of $\mathcal{V}$.\newline
If $\mathbf{x}$ and $\widetilde{\mathbf{x}}$ are two elements of $G\Omega
\left( \mathbb{R}^{d}\right) ^{p}$ such that $\left\| \mathbf{x}%
_{s,t}^{-1}\otimes \widetilde{\mathbf{x}}_{s,t}\right\| ^{p}\leq \varepsilon
\omega (s,t)$, then the corresponding solution of equation (\ref{DE}) $%
\mathbf{z}$ and $\widetilde{\mathbf{z}}$ satisfy $\left\| \mathbf{z}%
_{s,t}^{-1}\otimes \widetilde{\mathbf{z}}_{s,t}\right\| ^{p}\leq \delta
\left( \varepsilon \right) \omega (s,t)$ where $\delta $ is a continuous
function such that $\delta (0)=0$ (i.e. the map $\mathbf{x\rightarrow z}$ is
continuous, and hence the It\^{o} map $\mathbf{x\rightarrow y}$ is
continuous).

\begin{proof}
A simple translation of the first and second level estimates in \cite{Ly} or 
\cite{LQ} to our norm $|||.|||$ or $\left\Vert .\right\Vert $
\end{proof}
\end{theorem}

\bigskip

A simple corollary of it, observed in \cite{F}, is the continuity of the It%
\^{o} map in ``H\"{o}lder type''\ norm. The same simple argument gives the
continuity of the It\^{o} map in ``modulus topologies''. First, we let 
\begin{equation*}
\Xi _{p}=\left\{ \varphi :[0,1]\rightarrow \mathbb{R}^{+},\text{ with }%
\varphi (0)=0\text{ and }\varphi ^{p}\text{ is strictly increasing and convex%
}\right\} .
\end{equation*}
Such set is obviously not empty, $t\rightarrow t^{1/p}$ being one example of
an element of $\Xi _{p}$. Let us look at some more complicated one.

\begin{example}
\label{ex}Let $\alpha >0$. Then, the function $x\rightarrow \left( x(-\ln
x)^{a}\right) ^{p/2}$ is strictly increasing and convex in a neighborhood of 
$0$ (it can be checked by differentiating it twice). Let $\chi
_{a,p}=\inf_{x>0}\frac{d^{2}\left( x(-\ln x)^{a}\right) ^{p/2}}{dx^{2}}<0.$
Then define 
\begin{equation*}
\phi _{a,p}(x)=\left\{ 
\begin{array}{c}
\sqrt{x(-\ln x)^{a}}\text{ if }x\in \lbrack 0,\chi _{a,p}] \\ 
\left( \phi _{a,p}^{p}(\chi _{a,p})+\left( \phi _{a,p}^{p}\right) ^{^{\prime
}}(\chi _{a,p})\left( x-\chi _{a,p}\right) \right) ^{1/p}\text{ if }x\in
\lbrack \chi _{a,p},1]\text{ }%
\end{array}
\right.
\end{equation*}
In other words, $\phi _{a,p}^{p}(x)$ is the smallest convex function
dominating $\left( x(-\ln x)^{a}\right) ^{p/2}$. Remark that for $a>1,$%
\begin{equation*}
\lim_{x\rightarrow 0}\frac{\phi _{1,p}(x)}{\phi _{a,p}(x)}=0.
\end{equation*}
\end{example}

For a function $\varphi \in \Xi _{p},$ we define a distance between two $%
G\left( \mathbb{R}^{d}\right) $-valued paths $\mathbf{x}$ and $\widetilde{%
\mathbf{x}}$ 
\begin{equation*}
d_{\varphi }(\mathbf{x},\widetilde{\mathbf{x}})=\sup_{0\leq s<t\leq 1}\frac{%
\left\| \mathbf{x}_{s,t}^{-1}\otimes \widetilde{\mathbf{x}}_{s,t}\right\| }{%
\varphi (t-s)}.
\end{equation*}
For a single $G\left( \mathbb{R}^{d}\right) $-valued path $\mathbf{x}$, we
let 
\begin{equation*}
\left\| \mathbf{x}\right\| _{\varphi }=\sup_{0\leq s<t\leq 1}\frac{\left\| 
\mathbf{x}_{s,t}\right\| }{\varphi (t-s)}.
\end{equation*}
We also let $d_{\infty }(\mathbf{x},\widetilde{\mathbf{x}})=\sup_{0\leq
s<t\leq 1}\left\| \mathbf{x}_{s,t}^{-1}\otimes \widetilde{\mathbf{x}}%
_{s,t}\right\| $ and $\left\| \mathbf{x}\right\| _{\infty }=\sup_{0\leq
s<t\leq 1}\left\| \mathbf{x}_{s,t}\right\| .$

It is straightforward to check that $d_{\varphi }$ and $d_{\infty }$ are
distances on the space $C_{0}\left( G\left( \mathbb{R}^{d}\right) \right) .$

\begin{corollary}
\label{corcon}Let $\mathbf{x,}\widetilde{\mathbf{x}}\in G\Omega \left( 
\mathbb{R}^{d}\right) ^{p}$, $\epsilon >0$ and $\mathcal{V}$ be a linear map
from $\mathbb{R}^{d}$ into the $Lip[p+\varepsilon ,\mathbb{R}^{n}]$ vector
fields, and $\mathbf{z,}\widetilde{\mathbf{z}}\in G\Omega (\mathbb{R}%
^{d}\oplus \mathbb{R}^{n})^{p}$ the corresponding solution of equation (\ref%
{DE}). There exists a continuous function $\delta $ satisfying $\delta (0)=0$%
, such that 
\begin{equation*}
d_{\varphi }(\mathbf{x},\widetilde{\mathbf{x}})\leq \varepsilon \Rightarrow
d_{\varphi }(\mathbf{z},\widetilde{\mathbf{z}})\leq \delta (\varepsilon ).
\end{equation*}
\end{corollary}

\begin{remark}
The requirement that $\varphi \in \Xi _{p}$ to define $d_{\varphi }$ is only
for convenience (so that $(s,t)\rightarrow \varphi ^{p}(t-s)$ is a control).
Indeed, if $\widetilde{\varphi }$ is another increasing function such that $%
\widetilde{\varphi }$ is equivalent to $\varphi $ at $0$, then the
topologies on $G\left( \mathbb{R}^{d}\right) $-valued paths space induced by 
$d_{\varphi }$ and $d_{\widetilde{\varphi }}$ are identical.
\end{remark}

\subsection{The Translation Operator on Rough Path Space\label{trans_op}}

We define the translation operator, first introduced in a more general
situation in \cite{Ly}. Let $q$ be real such that $1/q+1/p>1$.

The following definition is motivated by replacing $x$ by $x+f$ in (\ref%
{signature}).

\begin{definition}
\label{def_tr}Let $\mathbf{x}\in G\Omega \left( \mathbb{R}^{d}\right) ^{p}$
and $\mathbf{f}$ be a $G\left( \mathbb{R}^{d}\right) $-valued path of finite 
$q$-variation. We let $f_{t}=\pi _{1}\left( \mathbf{f}_{t}\right) $ and $%
x_{t}=\pi _{1}\left( \mathbf{x}_{t}\right) .$ Then define $T_{\mathbf{f}}(%
\mathbf{x)}$ by 
\begin{equation*}
\pi _{1}\left( T_{\mathbf{f}}(\mathbf{x)}\right) _{s,t}=f_{s,t}+x_{s,t}
\end{equation*}
and 
\begin{eqnarray*}
\pi _{2}\left( T_{\mathbf{f}}(\mathbf{x)}_{s,t}\right) &=&\mathbf{\pi }%
_{2}\left( \mathbf{f}_{s,t}\right) +\pi _{2}\left( \mathbf{x}_{s,t}\right) +%
\frac{1}{2}\int_{s}^{t}f_{s,u}\otimes dx_{u}+\frac{1}{2}\int_{s}^{t}x_{s,u}%
\otimes df_{u} \\
&&-\frac{1}{2}\int_{s}^{t}df_{s,u}\otimes x_{u}-\frac{1}{2}%
\int_{s}^{t}dx_{s,u}\otimes f_{u},
\end{eqnarray*}
where the integrals are well defined Young integrals.
\end{definition}

Remark that it is easily checked that $T_{\mathbf{f}}(\mathbf{x)}_{s,t}=T_{%
\mathbf{f}}(\mathbf{x)}_{0,s}^{-1}\otimes T_{\mathbf{f}}(\mathbf{x)}_{0,t}.$

\begin{theorem}
\label{weakcontinuity}Let $\mathbf{x}$ be a $G\left( \mathbb{R}^{d}\right) $%
-valued path of finite $p$-variation controlled by $\varepsilon ^{p}\omega
(s,t)$, and $\mathbf{f}$ be a $G\left( \mathbb{R}^{d}\right) $-valued path
of finite $q$-variation controlled by $\varkappa (s,t)$. Assume moreover
that $\varkappa \left( s,t\right) ^{1/q}\leq C\omega \left( s,t\right)
^{1/p} $. Then for all $s<t$ and $\varepsilon <1$, 
\begin{equation*}
d\left( T_{\mathbf{f}}(\mathbf{x)}_{s,t},\mathbf{f}_{s,t}\right) \leq C\sqrt{%
\varepsilon }\omega \left( s,t\right) ^{1/p}.
\end{equation*}
\end{theorem}

\begin{proof}
\begin{eqnarray*}
\mathbf{f}_{s,t}^{-1}\otimes \widetilde{T}_{\mathbf{f}}(\mathbf{x)}_{s,t}
&=&\exp \left( x_{s,t}-\frac{1}{2}[f_{s,t},x_{s,t}]+\pi _{2}\left( \mathbf{x}%
_{s,t}\right) \right) \\
&&\otimes \exp \left( \frac{1}{2}\int_{s}^{t}f_{s,u}\otimes dx_{u}+\frac{1}{2%
}\int_{s}^{t}x_{s,u}\otimes df_{u}\right) \\
&&\otimes \exp \left( -\frac{1}{2}\int_{s}^{t}df_{s,u}\otimes x_{u}-\frac{1}{%
2}\int_{s}^{t}dx_{s,u}\otimes f_{u}\right) .
\end{eqnarray*}
Hence, by inequality (\ref{boundnorm}) and Young inequality \cite{Yo} (which
says that for all $s<t$, $\left| \int_{s}^{t}f_{s,u}\otimes dx_{u}\right|
\leq C\varkappa \left( s,t\right) ^{1/q}\omega \left( s,t\right) ^{1/p}$ and
similar inequalities for the other Young integrals), we get that 
\begin{eqnarray*}
\left\| \mathbf{f}_{s,t}^{-1}\otimes \widetilde{T}_{\mathbf{f}}(\mathbf{x)}%
_{s,t}\right\| &\leq &\varepsilon \omega ^{1/p}(s,t)+\sqrt{\varepsilon
^{2}\omega ^{2/p}(s,t)+C\varepsilon \varkappa \left( s,t\right) ^{1/q}\omega
\left( s,t\right) ^{1/p}} \\
&\leq &C\sqrt{\varepsilon }\omega \left( s,t\right) ^{1/p}.
\end{eqnarray*}
\end{proof}

\section{The Enhanced Brownian Motion}

In what follows we will lift Brownian motion as $(\mathbb{R}^{d},+)$-valued
to a $\left( G\left( \mathbb{R}^{d}\right) ,\otimes \right) $-valued
process. Early work by Gaveau \cite{Ga} and in particular the presentation
in \cite{Ma} use related algebraic ideas. See also \cite{Ne}.

\subsection{Two Classical Properties of the Brownian Motion}

Let $\left( C_{0}(\mathbb{R}^{d}),\mathcal{F},\left( \mathcal{F}_{t}\right)
_{t},\mathbb{P}\right) $ be the Wiener space. The evaluation operator $B$ is
then under $\mathbb{P}$ a Brownian motion starting at $0$. \newline
Applying Garsia, Rodemich and Rumsey inequality, it is not too difficult to
see \cite{SV,VaLN} that for all $s<t$, 
\begin{equation*}
\left\| B_{s,t}\right\| \leq \frac{4}{\sqrt{\alpha }}\int_{0}^{t-s}\sqrt{%
\frac{\log \left( 1+\frac{4F}{u^{2}}\right) }{u}}du
\end{equation*}
where $F$ is a $L^{1}$-random variable and a constant $\alpha >0$,
sufficiently small.

We denote by $W^{1,2}$ the Cameron-Martin space 
\begin{equation*}
\left\{ h:[0,1]\rightarrow \mathbb{R}^{d},\text{ }h(t)=\int_{0}^{t}h^{^{%
\prime }}(t)dt\text{ with }h^{^{\prime }}\in L^{2}\left( \left[ 0,1\right]
\right) \right\} .
\end{equation*}
Cameron-Martin theorem, e.g. \cite{KS,RY}, states that if $f\in W^{1,2}$ is $%
\left( \mathcal{F}_{t}\right) _{t}$-adapted, then the law of $\left(
B_{t}\right) _{0\leq t\leq 1}$ (i.e. the Wiener probability $\mathbb{P}$)
and the law of $\left( B_{t}+f(t)\right) _{0\leq t\leq 1}$ (that we will
denote $\mathbb{P}^{f}$) are equivalent.

We will now extend $B$ to a geometric $p$-rough path, show that our
``enhanced''\ Brownian motion has a similar modulus of continuity and extend
Cameron-Martin theorem to our enhanced Brownian motion.

\subsection{Their Extensions to The Enhanced Brownian Motion}

\subsubsection{Definition of the Enhanced Brownian Motion\label{def_EBM}}

The Enhanced Brownian Motion was first defined in \cite{Sip}. See also \cite%
{LQ}.

For $n\in \mathbb{N}$ and $x$ a $\mathbb{R}^{d}$-valued path, we denote by $%
x^{n}$ the path which agrees with $x$ at the points $\frac{k}{2^{n}}%
,k=0,\ldots ,2^{n}$ and which is linear in the intervals $\left[ \frac{k}{%
2^{n}},\frac{k+1}{2^{n}}\right] ,k=0,\ldots ,2^{n}-1$. As $x^{n}$ has finite 
$1$-variation, we can define $\mathbf{x}^{n}=S(x^{n}),$ that is the natural $%
G\left( \mathbb{R}^{d}\right) $-valued path lying above $x^{n}$.

We denote by $\Gamma $ the almost surely defined map 
\begin{equation*}
B\rightarrow \left( t\rightarrow \lim_{n\rightarrow \infty
}S(B^{n})_{t}\right) ,
\end{equation*}
and we let $\left( \mathbf{B}_{t}\right) _{0\leq t\leq 1}$ be the $G\left( 
\mathbb{R}^{d}\right) $-valued path $\Gamma (B).$ Note that \cite{LQ} 
\begin{equation*}
\mathbf{B}_{t}=_{a.s.}\exp (B_{t}+\frac{1}{2}\int_{s}^{t}\left( B_{u}\otimes
\circ dB_{u}-\circ dB_{u}\otimes B_{u}\right) ),
\end{equation*}
where we have used Stratonovich integration . We will call $\left( \mathbf{B}%
_{t}\right) _{t\geq 0}$ the enhanced Brownian motion.

\begin{remark}
Almost surely, $\Gamma \circ \pi _{1}\left( \mathbf{B}\right) =\mathbf{B}$
and $\pi _{1}\circ \Gamma \left( B\right) =B.$
\end{remark}

\begin{remark}
The inverse of $\exp $, denoted by $\log $, provides a global chart for $%
G\left( \mathbb{R}^{d}\right) $. Gaveau's diffusion, \cite{Ma}, is exactly
our EBM seen through this chart.
\end{remark}

\subsubsection{Modulus of Continuity for the Enhanced Brownian Motion}

We now cite Garsia, Rodemich and Rumsey inequality, but the function $f$
below is $G\left( \mathbb{R}^{d}\right) $-valued (while usually, $f$ \ takes
values in a normed vector space). Using the properties $\left\| .\right\| $%
on $G\left( \mathbb{R}^{d}\right) $ the proofs are identical.

\begin{theorem}
\label{GRR}Let $\Psi $ and $p$ be continuous strictly increasing functions
on $[0,\infty )$ with $p(0)=\Psi (0)=0$ and $\Psi (x)\rightarrow \infty $ as 
$x\rightarrow \infty $. Given $f\in C_{0}\left( [0,1],G\left( \mathbb{R}%
^{d}\right) \right) $, if 
\begin{equation}
\int_{0}^{1}\int_{0}^{1}\Psi \left( \frac{f(s)^{-1}\otimes f(t)}{p(\left|
t-s\right| )}\right) dsdt\leq F,  \label{GRR_inequ}
\end{equation}
then for $0\leq s<t\leq 1,$%
\begin{equation*}
\left\| f(s)^{-1}\otimes f(t)\right\| \leq 8\int_{0}^{t-s}\Psi ^{-1}\left( 
\frac{4D}{u^{2}}\right) dp(u).
\end{equation*}
\end{theorem}

Applying the Garsia, Rodemich and Rumsey inequality, we obtain a modulus of
continuity for the enhanced Brownian motion

\begin{theorem}
Define $\zeta (x)=\frac{1}{2\sqrt{2}}\int_{0}^{x}\sqrt{\frac{\log \left( 1+%
\frac{1}{u^{2}}\right) }{u}}du$. Then there exists a random variable $Z\geq
1 $ $a.s.$ and in $L^{1},$ and a constant $C>0$ such that for all $s<t$, 
\begin{equation}
\left\| \mathbf{B}_{s,t}\right\| \leq CZ^{1/4}\zeta \left( \frac{t-s}{\sqrt{Z%
}}\right) .  \label{modulus_good}
\end{equation}
\end{theorem}

\begin{proof}
As in the proof of L\'{e}vy's modulus of continuity in \cite{SV,VaLN} we use
the Garsia, Rodemich and Rumsey inequality with $f\left( t\right) =\mathbf{B}%
_{t},$ $p(x)=\sqrt{x}$ and $\Psi (x)=\exp \left( \alpha x^{2}\right) $. We
obtain 
\begin{equation*}
\left\Vert \mathbf{B}_{s,t}\right\Vert \leq C\int_{0}^{t-s}\sqrt{\frac{\log
\left( 1+\frac{4F}{u^{2}}\right) }{u}}du,
\end{equation*}%
where $F$ is the (now random) left hand side of (\ref{GRR_inequ}). Since $%
\left\Vert \mathbf{B}_{s,t}\right\Vert \overset{\text{law}}{=}\sqrt{t-s}%
\left\Vert \mathbf{B}_{0,1}\right\Vert $ the expectation of $F$ is estimated
by 
\begin{equation}
\mathbb{E}\left( \exp (\alpha \left\Vert \mathbf{B}_{0,1}\right\Vert
^{2}\right) .  \label{integrability}
\end{equation}%
We claim that this last expression is finite for a small enough $\alpha >0$.
Remark that $\pi _{1}\left( \mathbf{B}\right) ,$ resp. $\pi _{2}\left( 
\mathbf{B}\right) $ are some elements of the first (resp. second) Wiener-It%
\^{o} Chaos$.$ By general integrability properties of the Wiener-It\^{o}
chaos \cite[p 207]{RY}, there exists $\tilde{\alpha}>0$ 
\begin{eqnarray*}
\mathbb{E}\left( \exp (\tilde{\alpha}\left\vert \pi _{1}\left( \mathbf{B}%
_{0,1}\right) \right\vert _{\mathbb{R}^{d}}^{2})\right) &<&\infty , \\
\mathbb{E}\left( \exp (\tilde{\alpha}\left\vert \pi _{2}\left( \mathbf{B}%
_{0,1}\right) \right\vert _{\mathbb{R}^{d}\otimes \mathbb{R}^{d}})\right)
&<&\infty .
\end{eqnarray*}%
By inequality (\ref{boundnorm}), the finiteness of the expectation of $\exp
\left( \alpha \left\Vert \mathbf{B}_{0,1}\right\Vert ^{2}\right) $ is easily
obtain\footnote{%
One could prove this directly as we know the density of $\mathbf{B}_{0,1}$%
\cite{Le}.}. It remains just to define $Z=\max \left\{ 4F,1\right\} $ (we
want $Z\geq 1$ for technical convenience later on), and do a change of
variable to obtain inequality (\ref{modulus_good}).
\end{proof}

\begin{remark}
As $\zeta (x)\sim _{x\rightarrow 0}\sqrt{-x\ln x},$%
\begin{equation}
\overline{\lim_{\delta \rightarrow 0}}\sup_{\substack{ 0\leq s<t\leq 1  \\ %
\left| t-s\right| \leq \delta }}\frac{\left\| \mathbf{B}_{s,t}\right\| }{%
\sqrt{-\delta \ln \delta }}\leq C.  \label{modulus_limit}
\end{equation}
On the other hand we can trivially get a deterministic lower bound by noting 
$\left| B_{s,t}\right| _{\mathbb{R}^{d}}\leq \left\| \mathbf{B}%
_{s,t}\right\| $ and using L\'{e}vy's result. All this is known (with an
equality) for hypoelliptic diffusions on Nilpotent group in \cite{Ba1} and
elliptic diffusions in \cite{Ba2}, using a natural metric associated to the
diffusion.
\end{remark}

\begin{lemma}
\label{lem_ine_eta}There exists a constant $C$ such that for all $x,y\in
\lbrack 0,1]$, 
\begin{equation*}
\zeta (xy)\leq C\zeta (x)\zeta (y).
\end{equation*}
\end{lemma}

\begin{proof}
For $a$ small enough, there exists constants $K_{1}$ and $K_{2}$ such that
if $x\in (0,a]$%
\begin{equation*}
K_{1}\sqrt{-x\ln x}\leq \zeta (x)\leq K_{2}\sqrt{-x\ln x}.
\end{equation*}
Hence, from the inequality $\forall x,y\in (0,a],$ $-\ln (xy)\leq \frac{%
-2\ln a}{\ln ^{2}a}\ln (x)\ln (y),$ we obtain that for all $x,y\in \lbrack
0,a]$, $\zeta (xy)\leq C\zeta (x)\zeta (y)$ for a constant $C$. \newline
For a fixed $b$, it is easily seen that $0<\inf_{x\in (0,\frac{1}{b})}\frac{%
\zeta (xb)}{\zeta (x)}<\sup_{x\in (0,\frac{1}{b})}\frac{\zeta (xb)}{\zeta (x)%
}<\infty $ . Hence, if $x,y\in \lbrack 0,1],$%
\begin{eqnarray*}
\zeta (xy) &=&\zeta \left( \frac{axay}{a^{2}}\right) \leq \sup_{z\in
(0,a^{2})}\frac{\zeta (z/a^{2})}{\zeta (z)}\zeta \left( axay\right) \\
&\leq &C\sup_{z\in (0,a^{2})}\frac{\zeta (z/a^{2})}{\zeta (z)}\zeta \left(
ay\right) \zeta \left( ay\right) \\
&\leq &C\frac{\sup_{z\in (0,a^{2})}\frac{\zeta (z/a^{2})}{\zeta (z)}}{\left(
\inf_{z\in (0,1/a)}\frac{\zeta (za)}{\zeta (z)}\right) ^{2}}\zeta (x)\zeta
(y).
\end{eqnarray*}
\end{proof}

\begin{proposition}
\label{almost_good_control} 
\begin{equation}
\left\Vert \mathbf{B}_{s,t}\right\Vert \leq M\zeta \left( t-s\right) ,
\label{control_Bst_with_M}
\end{equation}
where $M$ is a random variable for which there exists a constant $\lambda >0$
such that $E(\exp (\lambda M^{2}))<\infty .$
\end{proposition}

\begin{proof}
From the previous lemma, we can set 
\begin{equation*}
M=CZ^{1/4}\zeta \left( \frac{1}{\sqrt{Z}}\right) =2C\int_{0}^{1}\sqrt{\ln
\left( 1+\frac{Z}{v^{4}}\right) }dv.
\end{equation*}
Hence, by Jensen inequality, 
\begin{eqnarray*}
E(\exp (\lambda M^{2})) &\leq &E\left( \int_{0}^{1}\exp \left( \left(
2\lambda C\right) ^{2}\ln \left( 1+\frac{Z}{v^{4}}\right) \right) dv\right)
\\
&\leq &E\left( \int_{0}^{1}\left( 1+\frac{Z}{v^{4}}\right) ^{\left( 2\lambda
C\right) ^{2}}dv\right) \\
&\leq &E\left( \int_{0}^{1}\left( \frac{2Z}{v^{4}}\right) ^{\left( 2\lambda
C\right) ^{2}}dv\right) \\
&<&\infty
\end{eqnarray*}
when $\lambda <\frac{1}{4C}$.
\end{proof}

The last estimate looks like a control but $\zeta ^{p\text{ }}$is not convex
on the entire interval $[0,1]$. We define $\phi _{p}$ to be $p^{th}$-root of
the smallest convex function dominating $x\rightarrow \zeta \left( x\right)
^{p}$ (remark that $\phi _{p}$ is very similar to $\phi _{1,p}$ of example %
\ref{ex}, as $\frac{\zeta \left( x\right) }{\sqrt{-x\ln x}}\rightarrow
_{x\rightarrow 0}1$).

\begin{corollary}
\label{good_control}The $p$-variation of $\mathbf{B}$ is controlled by $%
(s,t)\rightarrow C^{p}Z^{p/4}\phi _{p}^{p}\left( \frac{t-s}{\sqrt{Z}}\right) 
$ and also by $(s,t)\rightarrow M^{p}\phi _{p}^{p}\left( t-s\right) $.
\end{corollary}

\subsubsection{A Cameron-Martin Theorem on the Group}

This section, despite being short and quite trivial, will be crucial in the
proof of the support theorem. For $f\in W^{1,2}$, $T_{f}(\mathbf{B})=T_{%
\mathbf{f}}(\mathbf{B})$ is defined by extending canonically $f$ (which has
finite $1$-variation) to a $G\left( \mathbb{R}^{d}\right) $-valued path $%
\mathbf{f}$ (of finite $1$-variation).

\begin{theorem}
\label{CM}Let $f\in W^{1,2}$ be an $\left( \mathcal{F}_{t}\right) $-adapted
path. Then the law of $\mathbf{B}$ is equivalent to the law of $T_{\mathbf{f}%
}(\mathbf{B}).$
\end{theorem}

\begin{proof}
By the definition \ref{def_tr} and properties of Young and Stratonovich
integral, $T_{\mathbf{f}}(\mathbf{B})=\Gamma (f+B).$ Hence, the law of $%
\mathbf{B}$ is $\mathbb{P}\circ \Gamma ^{-1}$ while the law of $T_{\mathbf{f}%
}(\mathbf{B})$ is $\mathbb{P}^{f}\circ \Gamma ^{-1}.$ Hence, by the
Cameron-Martin theorem, these two laws are equivalent.
\end{proof}

\section{Modulus of Continuity for Solution of SDEs\label{MC}}

\begin{theorem}
Let $y_{t}$ be the solution of the Stratonovich differential equation 
\begin{equation*}
dy_{t}=f_{0}\left( t,y_{t}\right) dt+f(t,y_{t})\circ dB_{t},
\end{equation*}
where $f_{0},f$ are $2+\varepsilon $-H\"{o}lder. For $h\in W^{1,2}$, we
denote by $F(h)$ the solution of the ordinary differential equation\footnote{%
For $h\in W^{1,2}$ but not piecewise $C^{1}$ this ODE\ still makes sense as
\ rough differential equation with driving signal of finite $p=1$ variation,
see \cite{Lej}.} 
\begin{equation*}
dF(h)_{t}=f_{0}\left( t,F(h)_{t}\right) dt+f(t,F(h)_{t})dh_{t}.
\end{equation*}
We also denote by $\mathbf{F}$ the extension of the\ It\^{o} map $F$ to the
space of geometric $p$-rough path. There exists a constant $C$ such that 
\begin{equation*}
\overline{\lim_{h\rightarrow 0}}\sup_{\substack{ 0\leq s<t\leq 1  \\ \left|
t-s\right| \leq \delta }}\frac{\left\| \mathbf{F(B)}_{s,t}\right\| }{\sqrt{%
-\delta \ln \delta }}\leq C\text{ \ \ \ a.s.}
\end{equation*}
In particular, 
\begin{equation*}
\overline{\lim_{h\rightarrow 0}}\sup_{\substack{ 0\leq s<t\leq 1  \\ \left|
t-s\right| \leq \delta }}\frac{\left| y_{t}-y_{s}\right| _{\mathbb{R}^{d}}}{%
\sqrt{-\delta \ln \delta }}\leq C\text{ \ \ \ a.s.}
\end{equation*}
\end{theorem}

\begin{proof}
By theorem \ref{ItoMap} and corollary \ref{good_control}, for $s<t$ such
that $C^{p}Z^{\frac{p}{4}}\phi _{p}^{p}\left( \frac{t-s}{\sqrt{Z}}\right)
<1, $ 
\begin{equation*}
\left\| \mathbf{F(B)}_{s,t}\right\| \leq K_{\varepsilon ,V}CZ^{\frac{1}{4}%
}\phi _{p}\left( \frac{t-s}{\sqrt{Z}}\right) ,
\end{equation*}
and we obtain our theorem by remarking once again that $\phi _{p}\left(
x\right) \sim _{x\rightarrow 0}\sqrt{-x\ln (x)}.$ The second inequality is
obvious from the fact that $\pi _{1}\left( \mathbf{F(B)}_{s,t}\right)
=y_{t}-y_{s}$.
\end{proof}

\section{On the Support Theorem}

We are going to show a support theorem for the enhanced Brownian motion.
Using the continuity of the It\^{o} map, we will recover the classical
support theorem (and even more). First, we need to look carefully at the
convergence of piecewise linear approximation of our Brownian motion to the
Enhanced Brownian motion, in various topologies.

\subsection{Convergence of some Smooth \textbf{paths }to the Enhanced
Brownian Motion}

\begin{proposition}
\label{good_est}Let $\varphi ,\phi _{p}\in \Xi _{p},$ such that $%
\lim_{x\rightarrow 0}\frac{\phi _{p}(x)}{\varphi (x)}=0$. Let $\mathbf{x,y}$
be two $G\left( \mathbb{R}^{d}\right) $-valued paths. Then, for all $A\geq
d_{\phi _{p}}(\mathbf{x},\mathbf{y}),$%
\begin{equation*}
d_{\varphi }(\mathbf{x},\mathbf{y})\leq A\left( \frac{\varphi ^{-1}}{\phi
_{p}^{-1}}\left( \frac{d_{\infty }\left( \mathbf{x},\mathbf{y}\right) }{A}%
\right) \right) ^{1/p}
\end{equation*}
\end{proposition}

Remark that the functions $\varphi =\phi _{a,p}$, $a>1$ of example \ref{ex}
satisfy the condition $\lim_{x\rightarrow 0}\frac{\phi _{p}(x)}{\varphi (x)}%
=0$

\begin{proof}
For all $s<t,$ 
\begin{eqnarray*}
\varphi ^{-1}\left( \frac{\left\Vert \mathbf{x}_{s,t}^{-1}\otimes \mathbf{y}%
_{s,t}\right\Vert }{A}\right) &=&\frac{\varphi ^{-1}}{\phi _{p}^{-1}}\left( 
\frac{\left\Vert \mathbf{x}_{s,t}^{-1}\otimes \mathbf{y}_{s,t}\right\Vert }{A%
}\right) \phi _{p}^{-1}\left( \frac{\left\Vert \mathbf{x}_{s,t}^{-1}\otimes 
\mathbf{y}_{s,t}\right\Vert }{A}\right) \\
&\leq &\left( \frac{\varphi ^{-1}}{\phi _{p}^{-1}}\left( \frac{d_{\infty
}\left( \mathbf{x},\mathbf{y}\right) }{A}\right) \right) \phi
_{p}^{-1}\left( \frac{d_{\phi _{p}}(\mathbf{x,y})}{A}\phi _{p}\left(
t-s\right) \right) \\
&\leq &\left( \frac{\varphi ^{-1}}{\phi _{p}^{-1}}\left( \frac{d_{\infty
}\left( \mathbf{x},\mathbf{y}\right) }{A}\right) \right) \left( t-s\right)
\end{eqnarray*}%
Therefore, for all $s<t$, 
\begin{eqnarray*}
\frac{\left\Vert x_{s,t}^{-1}\otimes y_{s,t}\right\Vert }{A} &\leq &\varphi
\left( \left( \frac{\varphi ^{-1}}{\phi _{p}^{-1}}\left( \frac{d_{\infty
}\left( \mathbf{x},\mathbf{y}\right) }{A}\right) \right) \left( t-s\right)
\right) \\
&\leq &\left( \frac{\varphi ^{-1}}{\phi _{p}^{-1}}\left( \frac{d_{\infty
}\left( \mathbf{x},\mathbf{y}\right) }{A}\right) \right) ^{1/p}\varphi
\left( t-s\right)
\end{eqnarray*}%
by convexity of $\varphi ^{p}$. Hence,

\begin{equation*}
d_{\varphi }(\mathbf{x},\mathbf{y})\leq A\left( \frac{\varphi ^{-1}}{\phi
_{p}^{-1}}\left( \frac{d_{\infty }\left( \mathbf{x},\mathbf{y}\right) }{A}%
\right) \right) ^{1/p}.
\end{equation*}
\end{proof}

\begin{corollary}
\label{crucial_conv_cor}Let $\varphi \in \Xi _{p},$ such that $%
\lim_{x\rightarrow 0}\frac{\phi _{p}(x)}{\varphi (x)}=\lim_{x\rightarrow 0}%
\frac{\sqrt{-x\ln x}}{\varphi (x)}=0$. Let $\mathbf{x}_{n}$ be a sequence of 
$G\left( \mathbb{R}^{d}\right) $-valued paths, which converges pointwise to
another $G\left( \mathbb{R}^{d}\right) $-valued path $\mathbf{x}$. Assume
that $\sup_{n}\left\Vert \mathbf{x}_{n}\right\Vert _{\phi _{p}}<\infty $ .
Then $\mathbf{x}_{n}$ converges to $\mathbf{x}$ in the topology induced by $%
d_{\varphi }$.
\end{corollary}

\begin{proof}
First notice that by Arzela-Ascoli theorem, 
\begin{equation*}
\widetilde{d_{\infty }}(\mathbf{x}_{n},\mathbf{x}):=\sup_{t\in \lbrack
0,1]}\left\Vert \mathbf{x}_{n,t}^{-1}\otimes \mathbf{x}_{t}\right\Vert
\rightarrow _{n\rightarrow \infty }0.
\end{equation*}%
But inequality (\ref{useful_ineq}) gives 
\begin{equation}
\widetilde{d_{\infty }}(\mathbf{x}_{n},\mathbf{x})\leq d_{\infty }(\mathbf{x}%
_{n},\mathbf{x})\leq C\left( \widetilde{d_{\infty }}(\mathbf{x}_{n},\mathbf{x%
})+\sqrt{\widetilde{d_{\infty }}(\mathbf{x}_{n},\mathbf{x})\left\Vert 
\mathbf{x}\right\Vert _{\infty }}\right) ,  \label{twonotionsofsupdistance}
\end{equation}%
hence $d_{\infty }\left( \mathbf{x}_{n},\mathbf{x}\right) \rightarrow
_{n\rightarrow \infty }0$. Therefore, $d_{\varphi }(\mathbf{x}_{n},\mathbf{x}%
),$ being bounded by 
\begin{equation*}
2\sup_{n}\left\Vert \mathbf{x}_{n}\right\Vert _{\phi _{p}}\left( \frac{%
\varphi ^{-1}}{\phi _{p}^{-1}}\left( \frac{d_{\infty }\left( \mathbf{x}_{n},%
\mathbf{x}\right) }{2\sup_{n}\left\Vert \mathbf{x}_{n}\right\Vert _{\phi
_{p}}}\right) \right) ^{1/p}
\end{equation*}%
goes to $0$ when $n$ tends to infinity.
\end{proof}

This corollary is going to allow us to prove that various approximations of
enhanced Brownian motion converge in the topology induced by $d_{\varphi }$,
where $\varphi \in \Xi _{p}$ is such that $\phi _{p}(x)=o(\varphi (x))$ as $%
x\rightarrow 0$. To obtain an accurate uniform control of $\left\| \mathbf{B}%
^{n}\right\| _{\phi _{p}}$, we first need the following result, in the
spirit of Doob's martingale inequality.

\begin{lemma}
Let $X$ a random variable such that $E\left( \exp \lambda X^{2}\right)
<\infty $, and $\mathcal{G}_{n}$ a sequence of $\sigma $-algebras. Define $%
X_{n}^{2}=E\left( X^{2}/\mathcal{G}_{n}\right) $. Then 
\begin{equation*}
E\left( \sup_{n}\exp \lambda X_{n}^{2}\right) <\infty \text{.}
\end{equation*}
\end{lemma}

\begin{proof}
Using Fubini and Doob's $L^{p}$ inequality, we obtain 
\begin{eqnarray*}
E\left( \sup_{n}\exp \lambda X_{n}^{2}\right) &=&1+\sum_{k=1}^{\infty }\frac{%
\lambda ^{k}}{k!}E\left( \sup_{n}X_{n}^{2k}\right) \\
&\leq &1+\sum_{k=1}^{\infty }\frac{\lambda ^{k}}{k!}\left( \frac{2k}{2k-1}%
\right) ^{2k}E\left( X^{2k}\right) \\
&\leq &2eE\left( \exp \lambda X^{2}\right) .
\end{eqnarray*}
\end{proof}

\begin{proposition}
\label{uniform_control}$t\rightarrow \mathbf{B}_{t}$ and $t\rightarrow 
\mathbf{B}_{t}^{n}$ are almost surely of finite $p$-variation uniformly
controlled by $(s,t)\rightarrow K^{p}\phi _{p}^{p}(t-s),$ where $K$ is a
random variable such that for $\lambda $ small enough $E(\exp (\lambda
K^{2}))<\infty $.
\end{proposition}

\begin{proof}
$\mathbf{B}$ is controlled in $p$-variation by $C^{p}M^{p}\phi _{p}^{p}(t-s)$%
, where, for $\lambda $ small enough $E(\exp (\lambda M^{2}))<\infty $. By
inequality (\ref{boundnorm}), we see that for all $s<t$ 
\begin{equation*}
|||\mathbf{B}_{s,t}|||\leq \frac{M}{c_{1}}\phi _{p}\left( t-s\right) .
\end{equation*}
Now define $\mathcal{G}_{n}$ the $\sigma $-algebra generated by the random
variables $B_{\frac{k}{2^{n}}}.$ Then, $\log \mathbf{B}_{s,t}^{n}=\mathbb{E}%
\left( \log \mathbf{B}_{s,t}/\mathcal{G}_{n}\right) $ \cite{Ma}. Hence, 
\begin{eqnarray*}
|||\mathbf{B}_{s,t}^{n}||| &=&|||\exp \mathbb{E}\left( \log \mathbf{B}_{s,t}/%
\mathcal{G}_{n}\right) ||| \\
&=&\left| \mathbb{E}\left( \pi _{1}\left( \mathbf{B}_{s,t}\right) /\mathcal{G%
}_{n}\right) \right| _{\mathbb{R}^{d}}+\sqrt{\left| \mathbb{E}\left( \pi
_{2}\left( \mathbf{B}_{s,t}\right) /\mathcal{G}_{n}\right) \right| _{\mathbb{%
R}^{d}\otimes \mathbb{R}^{d}}} \\
&\leq &\sqrt{\mathbb{E}\left( \left| \pi _{1}\left( \mathbf{B}_{s,t}\right)
\right| _{\mathbb{R}^{d}}^{2}/\mathcal{G}_{n}\right) }+\sqrt{\mathbb{E}%
\left( \left| \pi _{2}\left( \mathbf{B}_{s,t}\right) \right| _{\mathbb{R}%
^{d}\otimes \mathbb{R}^{d}}/\mathcal{G}_{n}\right) } \\
&\leq &\sqrt{2}\sqrt{\mathbb{E}\left( |||\mathbf{B}_{s,t}|||^{2}/\mathcal{G}%
_{n}\right) } \\
&\leq &\frac{\sqrt{2}C}{c_{1}}\phi _{p}\left( t-s\right) \sqrt{\mathbb{E}%
\left( M^{2}/\mathcal{G}_{n}\right) }
\end{eqnarray*}
Define $\widetilde{K}^{2}=\sup_{n}\mathbb{E}\left( M^{2}/\mathcal{G}%
_{n}\right) $. For all $n$ and $s<t$, $|||\mathbf{B}_{s,t}^{n}|||\leq \frac{%
\sqrt{2}}{c_{1}}\widetilde{K}\phi _{p}\left( t-s\right) .$ By the previous
lemma, for $\lambda $ small enough $E(\exp (\lambda \widetilde{K}%
^{2}))<\infty $. A\ last use of (\ref{boundnorm}) gives us the proposition.
\end{proof}

As a consequence of proposition \ref{uniform_control} and corollary \ref%
{crucial_conv_cor}, we obtain the following corollaries.

\begin{corollary}
\label{conv_dphi}Let $\varphi \in \Xi _{p},$ such that $\lim_{x\rightarrow 0}%
\frac{\sqrt{-x\ln \left( x\right) }}{\varphi (x)}=0.$ Then, $\mathbf{B}^{n}$
converges almost surely to $\mathbf{B}$ in the topology induced by $%
d_{\varphi }.$
\end{corollary}

\begin{corollary}
\label{conv_dphi_adapted}Let $\varphi \in \Xi _{p},$ such that $%
\lim_{x\rightarrow 0}\frac{\sqrt{-x\ln \left( x\right) }}{\varphi (x)}=0.$
We define $B^{(n)}$ to be the $\left( \mathcal{F}_{t}\right) $-adapted path
such that $B_{t}^{(n)}=B_{\frac{\left[ 2^{n}t\right] }{2^{n}}}+\left( t-%
\frac{\left[ 2^{n}t\right] }{2^{n}}\right) \left( B_{\frac{\left[ 2^{n}t%
\right] }{2^{n}}}-B_{\frac{\left[ 2^{n}t\right] -1}{2^{n}}\vee 0}\right) .$
Then, $\mathbf{B}^{(n)}=S(B^{(n)})$ converges almost surely to $\mathbf{B}$
in the topology induced by $d_{\varphi }.$
\end{corollary}

\begin{proof}
We have seen in proposition \ref{uniform_control} that, almost surely, $%
\sup_{n}\left\| \mathbf{B}^{n}\right\| _{\phi _{p}}<\infty $ a.s. Observe
that $\mathbf{B}^{(n)}$ is essentially $\mathbf{B}^{n}$ shifted by $\epsilon
=1/2^{n}.$ More precisely, 
\begin{equation*}
\left\| \mathbf{B}_{s,t}^{(n)}\right\| =\left\| \mathbf{B}_{(s-\epsilon
)\wedge 0,(t-\epsilon )\wedge 0}^{n}\right\| \leq \sup_{n}\left\| \mathbf{B}%
^{n}\right\| _{\phi _{p}}\phi _{p}(t-s),
\end{equation*}

By corollary \ref{crucial_conv_cor} it suffices to show that, almost surely, 
$\mathbf{B}_{t}^{(n)}\rightarrow \mathbf{B}_{t}$ for fixed $t$. But this
simply follows from 
\begin{equation*}
\left\| \mathbf{B}_{t}^{-1}\otimes \mathbf{B}_{t}^{(n)}\right\| \leq 
\underset{\rightarrow 0}{\underbrace{\left\| \mathbf{B}_{t}^{-1}\otimes 
\mathbf{B}_{t}^{n}\right\| }}+\underset{\leq C\phi _{p}(\epsilon
)\rightarrow 0\text{ \ with }n\rightarrow \infty .}{\underbrace{\left\|
\left( \mathbf{B}_{t}^{n}\right) ^{-1}\otimes \mathbf{B}_{(t-\epsilon
)\wedge 0}^{n}\right\| }}.
\end{equation*}
\end{proof}

In particular, we recover the convergence in the $1/p$-H\"{o}lder distance
(and hence in the $p$-variation topology) of $\mathbf{B}^{n}$ and $\mathbf{B}%
^{(n)}$ to $\mathbf{B}$.

\subsection{Some more Convergence Results}

\begin{lemma}
\label{forsupport}Let $\varphi \in \Xi _{p},$ such that $\lim_{x\rightarrow
0}\frac{\sqrt{-x\ln \left( x\right) }}{\varphi (x)}=0$ and $f\in W^{1,2}$ be
an $\left( \mathcal{F}_{t}\right) $-adapted path. Then $\mathbb{P}$-almost
surely, $T_{f-B^{(n)}}(\mathbf{B})$ converges in the topology induced by $%
d_{\varphi }$ to $\mathbf{f}.$
\end{lemma}

\begin{proof}
Assume that we have shown that $T_{-B^{(n)}}(\mathbf{B})$ converges to $\exp
(0)$ in the topology induced by $d_{\varphi }$, i.e. that there exists a
sequence $\varepsilon _{n},$ which converges almost surely to $\exp (0),$
and such that $\varepsilon _{n}^{p}\varphi ^{p}(t-s)$ controls the $p$%
-variation of $T_{-B^{(n)}}(\mathbf{B})$. Then, as $\int_{s}^{t}\left|
f_{u}^{\prime }\right| du\leq \sqrt{t-s}\sqrt{\int_{0}^{1}\left|
f_{u}^{\prime }\right| ^{2}du}\leq C\sqrt{\int_{0}^{1}\left| f_{u}^{\prime
}\right| ^{2}du}\varphi (t-s),$ we obtain from theorem \ref{weakcontinuity}
and the equality $T_{f-B^{(n)}}(\mathbf{B})=T_{f}[T_{-B^{(n)}}(\mathbf{B})]$
that $d_{\varphi }\left( T_{f-B^{(n)}}(\mathbf{B}),\mathbf{f}\right) \leq C%
\sqrt{\varepsilon _{n}}$. Hence, we can assume that $f=0$.

Note that 
\begin{equation*}
T_{-B^{(n)}}(\mathbf{B})_{s,t}=\exp \left( B_{s,u}-B_{s,t}^{(n)}+\int_{s}^{t}%
\left[ B_{s,u}-B_{s,u}^{(n)},\circ d\left( B_{u}-B_{u}^{(n)}\right) \right]
\right)
\end{equation*}
Conditioning partially, i.e. only with respect to the values of the $i^{%
\text{th}}$ component of the Brownian motion at some fixed time, we obtain,
similarly as in the previous section (and as in \cite{F}), that, almost
surely, $T_{-B^{(n)}}(\mathbf{B})$ converges pointwise to $\exp (0)$, and
that for all $n$, 
\begin{equation*}
\left\| T_{-B^{(n)}}(\mathbf{B})_{s,t}\right\| \leq K^{\prime }\phi
_{p}(t-s),
\end{equation*}
where $K^{\prime }$ is a random variable such that for $\lambda $ small
enough, $E\left( \exp \lambda K^{\prime 2}\right) <\infty $. The proof is
then finish using corollary \ref{crucial_conv_cor}.
\end{proof}

\subsection{Support Theorem with Refined Norms}

As observed in \cite{F} we can combine the Rough Path approach of \cite{LQZ}
with ideas from \cite{Millet} and get an improved Support Theorem as a
corollary.

\begin{theorem}
The support of the law of $\mathbf{B}$ is the closure of $S(W^{1,2})$ in the
topology induced by $d_{\varphi }$, where $\varphi \in \Xi _{p}$ is such
that $\lim_{x\rightarrow 0}\frac{\sqrt{-x\ln \left( x\right) }}{\varphi (x)}%
=0.$
\end{theorem}

\begin{proof}
Corollary \ref{conv_dphi_adapted} implies classically that the support of
the law of $\mathbf{B}$ is contained in the closure of $S(W^{1,2})$ in the $%
d_{\varphi }$ topology. Reciprocally, our Cameron-Martin theorem \ref{CM}
implies it is enough to show that for a function $x\in W^{1,2}$, $%
T_{x-B^{(n)}}(\mathbf{B)}$ converges in the $d_{\varphi }$ topology to $S(x)$%
. But this was proven in lemma \ref{forsupport}.
\end{proof}

\begin{remark}
As noted in \cite{GNS}, the previous theorem would not work in the topology
induced by $d_{\phi _{p}}$. Indeed, $\left\{ \mathbf{x}\in C_{0}\left(
[0,1],G\left( \mathbb{R}^{d}\right) \right) ,\left\| \mathbf{x}\right\|
_{\phi _{p}}<\infty \right\} $ is not separable (but the set of continuous $%
G\left( \mathbb{R}^{d}\right) $-valued path $\mathbf{x}$ such that $\left\| 
\mathbf{x}\right\| _{\phi _{p}}<\infty $ and such that $\lim_{\delta
\rightarrow 0}\sup_{\left| t-s\right| \leq \delta }\frac{\left\| \mathbf{x}%
_{s,t}\right\| }{\varphi (\delta )}=0$ is separable).
\end{remark}

As in section \ref{MC}, we let $y_{t}$ be the solution of the Stratonovich
differential equation 
\begin{equation*}
dy_{t}=f_{0}\left( t,y_{t}\right) dt+f(t,y_{t})\circ dB_{t},
\end{equation*}
where $f_{0},f$ are $2+\varepsilon $-H\"{o}lder. For $h\in W^{1,2}$, we
denote by $F(h)$ the solution of the ordinary differential equation 
\begin{equation}
dF(h)_{t}=f_{0}\left( t,F(h)_{t}\right) dt+f(t,F(h)_{t})dh_{t}.  \label{ItoF}
\end{equation}
We also denote by $\mathbf{F}$ the extension of the\ It\^{o} map $F$ to the
space of geometric $p$-rough path. From the continuity of $\mathbf{F}$ in
the $d_{\varphi }$ topology (corollary \ref{corcon}), we instantaneously
obtain the following:

\begin{corollary}
The support of the law of $\Gamma (y)$ (the Stratonovich extension of $y$ to
a $p$-rough path) is the closure of $\mathbf{F}\left( S(W^{1,2})\right) $ in
the topology induced by $d_{\varphi }$, where $\varphi \in \Xi _{p}$ is such
that $\lim_{x\rightarrow 0}\frac{\sqrt{-x\ln \left( x\right) }}{\varphi (x)}%
=0.$
\end{corollary}

Projecting on the first level, we improve Stroock-Varadhan's result \cite%
{SV2}, its extension to H\"{o}lder norm \cite{Millet,BA1,BA2,ST} as well as
the $p$-variation result \cite{LQZ}. Our approach allows us to use more
refined topologies than the one induced by H\"{o}lder distances. Let $%
d_{1,\varphi }$ be the distance defined by the following formula: 
\begin{equation*}
d_{1,\varphi }(x,y)=\sup_{0\leq s<t\leq 1}\frac{\left|
y_{s,t}-x_{s,t}\right| _{\mathbb{R}^{d}}}{\varphi (t-s)}
\end{equation*}

\begin{corollary}
The support of the law of $y$ is the closure of $F(W^{1,2})$ in the topology
induced by $d_{1,\varphi }$, where $\varphi \in \Xi _{p}$ is such that $%
\lim_{x\rightarrow 0}\frac{\sqrt{-x\ln \left( x\right) }}{\varphi (x)}=0.$
\end{corollary}

\section{Large deviations results}

\subsection{Some preliminary results}

For $n\in \mathbb{N}$, we define the map $\Upsilon _{n}:C_{0}\left(
[0,1],G\left( \mathbb{R}^{d}\right) \right) \rightarrow C_{0}\left(
[0,1],G\left( \mathbb{R}^{d}\right) \right) $ where $\Upsilon _{n}\left( 
\mathbf{x}\right) $ is defined by:

$%
\begin{array}{ll}
(i) & \forall k\in \{0,\ldots ,2^{n}\},\text{ }\Upsilon _{n}\left( \mathbf{x}%
\right) _{\frac{k}{2^{n}}}=\mathbf{x}_{\frac{k}{2^{n}}} \\ 
(ii) & \forall k\in \{0,\ldots ,2^{n}\}\text{ and }\forall t\in \lbrack
0,2^{-n}]\text{, }\Upsilon _{n}\left( \mathbf{x}\right) _{\frac{k}{2^{n}},%
\frac{k}{2^{n}}+t}=\delta _{t2^{n}}\left( \mathbf{x}_{\frac{k}{2^{n}},\frac{%
k+1}{2^{n}}}\right) .%
\end{array}%
$

$\Upsilon _{n}\left( \mathbf{x}\right) $ is a piecewise linear approximation
of $\mathbf{x}$ to which we assign a non-canonical area!

\begin{lemma}
\label{uniform_bound_lemma}For all $n\in \mathbb{N}$, 
\begin{equation*}
\left\Vert \Upsilon _{n}\left( \mathbf{x}\right) \right\Vert _{\phi
_{p}}\leq C\left\Vert \mathbf{x}\right\Vert _{\phi _{p}}.
\end{equation*}
\end{lemma}

\begin{proof}
One can show, with similar techniques than in lemma \ref{lem_ine_eta}, that
for all $\alpha ,y\in (0,1],$%
\begin{equation}
\phi _{p}(\alpha y)\geq C\sqrt{\alpha }\phi _{p}(y).
\label{inequalitywithoutaname}
\end{equation}%
Whenever $0\leq s\leq t\leq 2^{-n},$ $\Upsilon _{n}\left( \mathbf{x}\right)
_{\frac{k}{2^{n}}+s,\frac{k}{2^{n}}+t}$ is equal to 
\begin{equation*}
\exp \left( 2^{n}(t-s)\pi _{1}\left( \mathbf{x}_{\frac{k}{2^{n}},\frac{k+1}{%
2^{n}}}\right) +2^{2n}\left( t^{2}-s^{2}\right) \pi _{2}\left( \mathbf{x}_{%
\frac{k}{2^{n}},\frac{k+1}{2^{n}}}\right) \right) .
\end{equation*}%
Hence, by inequality (\ref{boundnorm}),%
\begin{eqnarray}
\left\Vert \Upsilon _{n}\left( \mathbf{x}\right) _{\frac{k}{2^{n}}+s,\frac{k%
}{2^{n}}+t}\right\Vert &\leq &C\left\Vert \mathbf{x}_{\frac{k}{2^{n}},\frac{%
k+1}{2^{n}}}\right\Vert \left( 2^{n}(t-s)+2^{n/2}\sqrt{t-s}\right)
\label{ex_helix_control} \\
&\leq &C\left\Vert \mathbf{x}_{\frac{k}{2^{n}},\frac{k+1}{2^{n}}}\right\Vert
2^{n/2}\sqrt{t-s}  \notag
\end{eqnarray}%
Hence, for $\frac{k}{2^{n}}\leq s\leq t\leq \frac{k+1}{2^{n}},$%
\begin{eqnarray*}
\left\Vert \Upsilon _{n}\left( \mathbf{x}\right) _{s,t}\right\Vert &\leq
&C\left\Vert \mathbf{x}\right\Vert _{\phi _{p}}\phi _{p}\left( t-s\right) 
\frac{\phi _{p}\left( 2^{-n}\right) }{\phi _{p}\left( t-s\right) }\sqrt{%
\frac{t-s}{2^{-n}}} \\
&\leq &C\left\Vert \mathbf{x}\right\Vert _{\phi _{p}}\phi _{p}\left(
t-s\right) \text{,}
\end{eqnarray*}%
applying inequality (\ref{inequalitywithoutaname}) with $\alpha =\frac{t-s}{%
2^{-n}}$ and $y=2^{-n}$.\newline
For general $s\leq \frac{j}{2^{n}}\leq \frac{k}{2^{n}}\leq t$, as 
\begin{eqnarray}
\Upsilon _{n}\left( \mathbf{x}\right) _{s,t} &=&\Upsilon _{n}\left( \mathbf{x%
}\right) _{s,\frac{j}{2^{n}}}\otimes \Upsilon _{n}\left( \mathbf{x}\right) _{%
\frac{j}{2^{n}},\frac{k}{2^{n}}}\otimes \Upsilon _{n}\left( \mathbf{x}%
\right) _{\frac{k}{2^{n}},t}  \notag \\
&=&\Upsilon _{n}\left( \mathbf{x}\right) _{s,\frac{j}{2^{n}}}\otimes \mathbf{%
x}_{\frac{j}{2^{n}},\frac{k}{2^{n}}}\otimes \Upsilon _{n}\left( \mathbf{x}%
\right) _{\frac{k}{2^{n}},t}  \label{decompostiongamman}
\end{eqnarray}%
\begin{eqnarray*}
\left\Vert \Upsilon _{n}\left( \mathbf{x}\right) _{s,t}\right\Vert &\leq
&C\left\Vert \mathbf{x}\right\Vert _{\phi _{p}}\left( \phi _{p}\left( \frac{j%
}{2^{n}}-s\right) +\phi _{p}\left( \frac{k-j}{2^{n}}\right) +\phi _{p}\left(
t-\frac{k}{2^{n}}\right) \right) \\
&\leq &C\left\Vert \mathbf{x}\right\Vert _{\phi _{p}}\phi _{p}\left(
t-s\right) .
\end{eqnarray*}
\end{proof}

\begin{lemma}
\label{deterministic_convergence_lemma}For all $n\in \mathbb{N}$, 
\begin{equation*}
\frac{d_{\infty }\left( \mathbf{x},\Upsilon _{n}\left( \mathbf{x}\right)
\right) }{\left\Vert \mathbf{x}\right\Vert _{\phi _{p}}}\leq C\sqrt{\phi
_{p}\left( 2^{-n}\right) }
\end{equation*}
\end{lemma}

\begin{proof}
For $s\leq \frac{j}{2^{n}}\leq \frac{k}{2^{n}}\leq t$, we obtain using
equation (\ref{decompostiongamman}) 
\begin{equation*}
d\left( \mathbf{x}_{s,t},\Upsilon _{n}\left( \mathbf{x}\right) _{s,t}\right)
=d\left( \Upsilon _{n}\left( \mathbf{x}\right) _{s,\frac{j}{2^{n}}}\otimes 
\mathbf{x}_{\frac{j}{2^{n}},\frac{k}{2^{n}}}\otimes \Upsilon _{n}\left( 
\mathbf{x}\right) _{\frac{k}{2^{n}},t},\mathbf{x}_{s,\frac{j}{2^{n}}}\otimes 
\mathbf{x}_{\frac{j}{2^{n}},\frac{k}{2^{n}}}\otimes \mathbf{x}_{\frac{k}{%
2^{n}},t}\right) .
\end{equation*}%
We then use inequality (\ref{useful_ineq2}):%
\begin{eqnarray*}
d\left( \mathbf{x}_{s,t},\Upsilon _{n}\left( \mathbf{x}\right) _{s,t}\right)
&\leq &d\left( \Upsilon _{n}\left( \mathbf{x}\right) _{s,\frac{j}{2^{n}}},%
\mathbf{x}_{s,\frac{j}{2^{n}}}\right) +d\left( \Upsilon _{n}\left( \mathbf{x}%
\right) _{\frac{k}{2^{n}},t},\mathbf{x}_{\frac{k}{2^{n}},t}\right) \\
&&+\sqrt{d\left( \Upsilon _{n}\left( \mathbf{x}\right) _{s,\frac{j}{2^{n}}},%
\mathbf{x}_{s,\frac{j}{2^{n}}}\right) \left\Vert \mathbf{x}_{\frac{j}{2^{n}}%
,t}\right\Vert }.
\end{eqnarray*}%
Then we simply bound $d\left( \Upsilon _{n}\left( \mathbf{x}\right) _{s,%
\frac{j}{2^{n}}},\mathbf{x}_{s,\frac{j}{2^{n}}}\right) $ by $\left\Vert
\Upsilon _{n}\left( \mathbf{x}\right) _{s,\frac{j}{2^{n}}}\right\Vert
+\left\Vert \mathbf{x}_{s,\frac{j}{2^{n}}}\right\Vert \leq C\left\Vert 
\mathbf{x}\right\Vert _{\phi _{p}}\phi _{p}\left( 2^{-n}\right) $, and
similarly for $d\left( \Upsilon _{n}\left( \mathbf{x}\right) _{\frac{k}{2^{n}%
},t},\mathbf{x}_{\frac{k}{2^{n}},t}\right) $. Hence,%
\begin{equation*}
\frac{d\left( \mathbf{x}_{s,t},\Upsilon _{n}\left( \mathbf{x}\right)
_{s,t}\right) }{\left\Vert \mathbf{x}\right\Vert _{\phi _{p}}}\leq C\phi
_{p}\left( 2^{-n}\right) +C\sqrt{\phi _{p}\left( 2^{-n}\right) \phi
_{p}\left( 1\right) }.
\end{equation*}
\end{proof}

\begin{corollary}
\label{forldp}For all $n\in \mathbb{N}$, 
\begin{equation*}
d_{\varphi }(\mathbf{x},\Upsilon _{n}\left( \mathbf{x}\right) )\leq
C\left\Vert \mathbf{x}\right\Vert _{\phi _{p}}\left( \sup_{0\leq x\leq
C^{\prime }\sqrt{\phi _{p}\left( 2^{-n}\right) }}\frac{\varphi ^{-1}}{\phi
_{p}^{-1}}\left( x\right) \right) ^{1/p}.
\end{equation*}%
In particular, if $\left\Vert \mathbf{x}\right\Vert _{\phi _{p}}<\infty $, $%
d_{\varphi }(\mathbf{x},\Upsilon _{n}\left( \mathbf{x}\right) )$ converges
to $0$ when $n$ tends to infinity.
\end{corollary}

\begin{proof}
We apply proposition \ref{good_est}, with $A=C\left\Vert \mathbf{x}%
\right\Vert _{\phi _{p}}\leq d_{\phi _{p}}(\Upsilon _{n}\left( \mathbf{x}%
\right) ,\mathbf{x})$. It gives%
\begin{equation*}
d_{\varphi }(\mathbf{x},\Upsilon _{n}\left( \mathbf{x}\right) )\leq
C\left\Vert \mathbf{x}\right\Vert _{\phi _{p}}\left( \frac{\varphi ^{-1}}{%
\phi _{p}^{-1}}\left( \frac{d_{\infty }\left( \mathbf{x},\Upsilon _{n}\left( 
\mathbf{x}\right) \right) }{C\left\Vert \mathbf{x}\right\Vert _{\phi _{p}}}%
\right) \right) ^{1/p}.
\end{equation*}%
The result is then given by lemma \ref{deterministic_convergence_lemma}.
\end{proof}

\begin{remark}
A feature of this approximation is that it does not rely on dyadic (or
nested) approximations which are fundamental for our earlier martingale
approach. Indeed, the $2^{-n}$ appearing on the right hand side of above
estimate is readily replaced by the mesh of any dissection upon which $%
\Upsilon _{n}\left( \mathbf{x}\right) $ is constructed. Introducing the
right area in our approximation improves its convergence properties.
\end{remark}

\subsection{Schilder and Freidlin-Wentzell theorem with Refined Norms}

We extend Schilder theorem \cite{DS,DZ,BBK} to the Enhanced Brownian Motion
in our refined topology. First, we need the following lemma

\begin{lemma}
\label{continuity_LD}Let $\varphi \in \Xi _{p},$ such that $%
\lim_{x\rightarrow 0}\frac{\sqrt{x}}{\varphi (x)}=0.$ The maps 
\begin{equation*}
\Upsilon _{n}:\left( C_{0}\left( [0,1],G\left( \mathbb{R}^{d}\right) \right)
,d_{\infty }\right) \rightarrow \left( C_{0}\left( [0,1],G\left( \mathbb{R}%
^{d}\right) \right) ,d_{\varphi }\right)
\end{equation*}%
are continuous.
\end{lemma}

\begin{proof}
The map $\Upsilon _{n}$ from $\left( C_{0}\left( [0,1],G\left( \mathbb{R}%
^{d}\right) \right) ,d_{\infty }\right) $ into $\left( C_{0}\left(
[0,1],G\left( \mathbb{R}^{d}\right) \right) ,d_{\infty }\right) $ is clearly
continuous as easily seen using (\ref{useful_ineq2}) as before. Let $\mathbf{%
x}\in C_{0}\left( [0,1],G\left( \mathbb{R}^{d}\right) \right) $ and $s\leq 
\frac{j}{2^{n}}\leq \frac{k}{2^{n}}\leq t,$%
\begin{eqnarray*}
\left\Vert \Upsilon _{n}\left( \mathbf{x}\right) _{s,t}\right\Vert &\leq
&\left\Vert \Upsilon _{n}\left( \mathbf{x}\right) _{s,\frac{j}{2^{n}}%
}\right\Vert +\left\Vert \mathbf{x}_{\frac{j}{2^{n}},\frac{k}{2^{n}}%
}\right\Vert +\left\Vert \Upsilon _{n}\left( \mathbf{x}\right) _{\frac{k}{%
2^{n}},t}\right\Vert \\
&\leq &C\left\Vert \mathbf{x}_{\frac{j-1}{2^{n}},\frac{j}{2^{n}}}\right\Vert
2^{n/2}\sqrt{\frac{j}{2^{n}}-s}+C\left\Vert \mathbf{x}_{\frac{k-1}{2^{n}},%
\frac{k}{2^{n}}}\right\Vert 2^{n/2}\sqrt{t-\frac{k}{2^{n}}} \\
&&+1_{j<k}2^{n/2}\left\Vert \mathbf{x}_{\frac{j}{2^{n}},\frac{k}{2^{n}}%
}\right\Vert \sqrt{\frac{k}{2^{n}}-\frac{j}{2^{n}}} \\
&\leq &C2^{n/2}\left\Vert \mathbf{x}\right\Vert _{\infty }\sqrt{t-s}.
\end{eqnarray*}%
When $\frac{j}{2^{n}}\leq s\leq t\leq \frac{j+1}{2^{n}}$, we also have 
\begin{equation}
\left\Vert \Upsilon _{n}\left( \mathbf{x}\right) _{s,t}\right\Vert \leq
C\left\Vert \mathbf{x}_{\frac{k}{2^{n}},\frac{k+1}{2^{n}}}\right\Vert 2^{n/2}%
\sqrt{t-s},  \notag
\end{equation}%
as already notice in equation \ref{ex_helix_control}. Hence, 
\begin{equation*}
\sup_{s<t}\frac{\left\Vert \Upsilon _{n}\left( \mathbf{x}\right)
_{s,t}\right\Vert }{\sqrt{t-s}}\leq C_{n}\left\Vert \mathbf{x}\right\Vert
_{\infty }.
\end{equation*}%
The proof is then finished applying a slight modification (replacing $\phi
_{p}$ by $\sqrt{.}$) of proposition \ref{good_est}.
\end{proof}

\begin{theorem}
Let $\varphi \in \Xi _{p},$ such that $\lim_{x\rightarrow 0}\frac{\sqrt{%
-x\ln \left( x\right) }}{\varphi (x)}=0.$ The random variables $\delta
_{\varepsilon }\mathbf{B}$ satisfies a large deviation principle in the
topology induced by $d_{\varphi }$ with good rate function 
\begin{equation*}
I(\mathbf{x)=}\left\{ 
\begin{array}{l}
\frac{1}{2}\int_{0}^{1}\left| x_{u}^{\prime }\right| ^{2}du\text{, if }S(x)=%
\mathbf{x}\text{ for some }x\in W^{1,2} \\ 
+\infty \text{ otherwise.}%
\end{array}
\right.
\end{equation*}
\end{theorem}

\begin{proof}
The large deviation result in \cite{LQZ} tells us that $\delta _{\varepsilon
}\mathbf{B}$ satisfies a large deviation principle with good rate function $%
I $ using the topology induced by $d_{\infty }$. By the lemma \ref%
{continuity_LD}, if we prove that $\Upsilon _{n}\left( \delta _{\varepsilon }%
\mathbf{B}\right) $ is an exponentially good approximation of $\delta
_{\varepsilon }\mathbf{B}$, i.e. 
\begin{equation}
\lim_{n\rightarrow \infty }\overline{\lim_{\varepsilon \rightarrow 0}}%
\varepsilon ^{2}\log \mathbb{P}\left( d_{\varphi }\left( \Upsilon _{n}\left(
\delta _{\varepsilon }\mathbf{B}\right) ,\delta _{\varepsilon }\mathbf{B}%
\right) >\delta \right) =-\infty \text{.}  \label{exp_convergence}
\end{equation}%
and that for all $\alpha $, 
\begin{equation}
\lim_{n\rightarrow \infty }\sup_{\mathbf{x,}I(\mathbf{x)}\leq \alpha
}d_{\varphi }(\Upsilon _{n}\left( \mathbf{x}\right) ,\mathbf{x})=0\text{,}
\label{cond_DZ}
\end{equation}%
we will have shown our theorem, by applying theorem 4.2.23 in \cite{DZ}. 
\newline

Let us first prove equation (\ref{cond_DZ}). First observe that if $I(%
\mathbf{x)}\leq \alpha ,$ letting $x\in W^{1,2}$ be such that $S(x)=\mathbf{%
x,}$ we have, by definition of our homogeneous norm and by Cauchy-Schwartz, $%
\left\Vert \mathbf{x}_{s,t}\right\Vert \leq \int_{s}^{t}\left\vert
x_{u}^{\prime }\right\vert du\leq 2\sqrt{t-s}I(\mathbf{x)}$. Hence, $%
\left\Vert \mathbf{x}\right\Vert _{\phi _{p}}\leq 2\sup_{0\leq s<t\leq 1}%
\frac{\sqrt{t-s}}{\phi _{p}(t-s)}\sqrt{\alpha }=C\sqrt{\alpha }$. and hence
by corollary \ref{forldp}, 
\begin{equation*}
\sup_{\mathbf{x,}I(\mathbf{x)}\leq \alpha }d_{\varphi }(\Upsilon _{n}\left( 
\mathbf{x}\right) ,\mathbf{x})\leq C\sqrt{\alpha }\left( \sup_{0\leq x\leq
C^{\prime }\sqrt{\phi _{p}\left( 2^{-n}\right) }}\frac{\varphi ^{-1}}{\phi
_{p}^{-1}}\left( x\right) \right) ^{1/p}\rightarrow _{n\rightarrow \infty }0%
\text{.}
\end{equation*}%
\qquad

To obtain inequality (\ref{exp_convergence}), we just need to apply
corollary \ref{forldp} to the enhanced Brownian Motion. Indeed, 
\begin{eqnarray*}
\mathbb{P}\left( d_{\varphi }\left( \delta _{\varepsilon }\Upsilon
_{n}\left( \mathbf{B}\right) ,\delta _{\varepsilon }\mathbf{B}\right)
>\delta \right) &=&\mathbb{P}\left( d_{\varphi }\left( \Upsilon _{n}\left( 
\mathbf{B}\right) ,\mathbf{B}\right) >\frac{\delta }{\varepsilon }\right) \\
&\leq &\mathbb{P}\left( C\left\Vert \mathbf{B}\right\Vert _{\phi _{p}}\left(
\sup_{0\leq x\leq C^{\prime }\sqrt{\phi _{p}\left( 2^{-n}\right) }}\frac{%
\varphi ^{-1}}{\phi _{p}^{-1}}\left( x\right) \right) ^{1/p}>\frac{\delta }{%
\varepsilon }\right) \\
&\leq &\mathbb{P}\left( M>\frac{\delta }{\alpha _{n}\varepsilon }\right)
\end{eqnarray*}%
where $M$ is the constant in proposition \ref{almost_good_control}, and $%
\alpha _{n}=C\left( \sup_{0\leq x\leq C^{\prime }\sqrt{\phi _{p}\left(
2^{-n}\right) }}\frac{\varphi ^{-1}}{\phi _{p}^{-1}}\left( x\right) \right)
^{1/p}$ is a deterministic sequence which converges to $0$ when $%
n\rightarrow \infty $. As 
\begin{equation*}
\mathbb{P}\left( M>\frac{\delta }{\alpha _{n}\varepsilon }\right) \leq
E\left( \exp \left( \lambda M^{2}\right) \right) \exp \left( -\lambda \left( 
\frac{\delta }{\alpha _{n}\varepsilon }\right) ^{2}\right) ,
\end{equation*}%
\begin{eqnarray*}
\overline{\lim_{\varepsilon \rightarrow 0}}\varepsilon ^{2}\log \mathbb{P}%
\left( d_{\varphi }\left( \delta _{\varepsilon }\Upsilon _{n}\left( \mathbf{B%
}\right) ,\delta _{\varepsilon }\mathbf{B}\right) >\delta \right) &\leq &%
\overline{\lim_{\varepsilon \rightarrow 0}}\varepsilon ^{2}\log E\left( \exp
\left( \lambda M^{2}\right) \right) -\lambda \left( \frac{\delta }{\alpha
_{n}}\right) ^{2} \\
&\leq &-\lambda \left( \frac{\delta }{\alpha _{n}}\right) ^{2},
\end{eqnarray*}%
which gives inequality (\ref{exp_convergence}).
\end{proof}

By the continuity of the It\^{o} map in the topology induced by $d_{\varphi
} $, we obtain the following extension of Freidlin-Wentzell theorem \cite%
{DS,DZ,BL}.

\begin{corollary}
Let $\varphi \in \Xi _{p},$ such that $\lim_{x\rightarrow 0}\frac{\sqrt{%
-x\ln \left( x\right) }}{\varphi (x)}=0.$ Let $y_{t}^{\varepsilon }$ be the
solution of the Stratonovich differential equation 
\begin{equation*}
dy_{t}^{\varepsilon }=f_{0}\left( t,y_{t}^{\varepsilon }\right)
dt+\varepsilon f(t,y_{t}^{\varepsilon })\circ dB_{t},
\end{equation*}%
where $f_{0},f$ are $2+\alpha $-H\"{o}lder. The Stratonovich extension of $%
y^{\varepsilon }$ to a geometric rough path, i.e. $\Gamma
(y_{t}^{\varepsilon })=\mathbf{F}(\delta _{\varepsilon }\mathbf{B})$ ($%
\mathbf{F}$ has been defined in equation (\ref{ItoF})) satisfies a large
deviation principle in the topology induced by $d_{\varphi }$ with good rate
function 
\begin{equation*}
J(\mathbf{x)}=\inf_{\mathbf{F(y)=x}}I(\mathbf{y})\text{.}
\end{equation*}
\end{corollary}

Remark that if we only consider the first level of our paths $y^{\varepsilon
}$, we obtain the classical Freidlin-Wentzell theorem in the topology
induced by $d_{1,\varphi }$.

\subsection{Strassen Law}

A\ classical corollary of Schilder is the law of the iterated logarithm \cite%
{DS,DZ,LQZ}:

\begin{corollary}
Let $K=\left\{ \mathbf{x=}S(x),\text{ }x\in W^{1,2}\text{ and }%
\int_{0}^{1}\left\vert x_{u}^{\prime }\right\vert ^{2}du\leq 1\right\} $ and 
\begin{equation*}
\mathbf{\xi }_{t}^{n}=\delta _{\left( 2n\log \log n\right) ^{-1/2}}\mathbf{B}%
_{nt}\text{.}
\end{equation*}%
Then, if $\varphi \in \Xi _{p}$ is such that $\lim_{x\rightarrow 0}\frac{%
\sqrt{-x\ln \left( x\right) }}{\varphi (x)}=0$, 
\begin{equation*}
\lim_{n\rightarrow \infty }d_{\varphi }\left( \mathbf{\xi }^{n},K\right) =0,
\end{equation*}%
and the set of limit points of in $C_{0}\left( [0,1],G\left( \mathbb{R}%
^{d}\right) \right) $ with the topology induced by $d_{\varphi }$ is equal
to $K$.
\end{corollary}


\begin{thebibliography}{99}
\bibitem{Ba1} Baldi, P.; Chaleyat-Maurel, M.: Sur l'\'{e}quivalent du module
de continuit\'{e} des processus de diffusion. S\'{e}minaire de Probabilit%
\'{e}s, XXI, 404--427, Lecture Notes in Math., 1247, Springer, Berlin, 1987.

\bibitem{Ba2} Baldi, P.; Sanz-Sol\'{e}, M.: Modulus of continuity for
stochastic flows. Barcelona Seminar on Stochastic Analysis (St. Feliu de Gu%
\'{\i}xols, 1991), 1--20, Progr. Probab., 32, Birkh\"{a}user, Basel, 1993.

\bibitem{BBK} Baldi, P.; Ben Arous, G.; Kerkyacharian, G. Large deviations
and the Strassen theorem in H\"{o}lder norm. Stochastic Process. Appl. Vol
42, no. 1, 171--180. 1992.

\bibitem{BL} Ben Arous, G.; Ledoux, M. Grandes d\'{e}viations de
Freidlin-Wentzell en norme h\"{o}lderienne. S\'{e}minaire de Probabilit\'{e}%
s, XXVIII, 293--299, Lecture Notes in Math., 1583, Springer, Berlin, 1994.

\bibitem{BA1} Ben Arous, G\'{e}rard; Gruadinaru, Mihai Normes h\"{o}%
lderiennes et support des diffusions. C. R. Acad. Sci. Paris S\'{e}r. I
Math. 316, no. 3, 283--286. 1993.

\bibitem{BA2} Ben Arous, G\'{e}rard; Gruadinaru, Mihai; Ledoux, Michel H\"{o}%
lder norms and the support theorem for diffusions. Ann. Inst. H. Poincar\'{e}
Probab. Statist. 30, no. 3, 415--436. 1994.

\bibitem{Ch} Chen K.T.: Integration of Paths, geometric invariant, and a
generalized Campbell-Hausdorff formula. In : The Annals of Mathematics,
Second Series, Volume 65, Issue 1 : 163-178 (Jan 1957).

\bibitem{DS} Deuschel, J.-D., Stroock, D., : Large Deviations. Academic
Press, New York, 1989.

\bibitem{DZ} Dembo, Amir; Zeitouni, Ofer Large deviations techniques and
applications. Jones and Bartlett Publishers, Boston, MA, 1993.

\bibitem{ENO} Eddahbi, M.; N'zi, M; Ouknine, Y. : Grandes d\'{e}viations des
diffusions sur les espaces de Besov-Orlicz et application. Stochastics
Stochastics Rep. 65, no. 3-4, 299--315. 1999.

\bibitem{FS} Folland, G. B.; Stein, E. M. : Hardy spaces on homogeneous
groups. Princeton University Press, 1982.

\bibitem{F} Friz, P.K., Continuity of the It\^{o}-map for H\"{o}lder rough
path with applications to the Support Theorem in H\"{o}lder norm\ (preprint
2003).

\bibitem{Ga} Gaveau, B.: Principe de moindre action, propagation de la
chaleur et estim\'{e}es sous-elliptiques sur certains groupes nilpotents.
Acta Math. 139 (1-2), 95-153, 1977.

\bibitem{Gr} Gromov, M.: Carnot-Carath\'{e}odory spaces seen from within.
Sub-Riemannian geometry, 79-323, Progr. Math., 144, Birkhauser, Basel, 1996.

\bibitem{GNS} Gy\"{o}ngy, I.; Nualart, D.; Sanz-Sol\'{e}, M. : Approximation
and support theorems in modulus spaces. Probab. Theory Related Fields 101,
no. 4, 495--509. 1995.

\bibitem{HL} Hambly, B. M.; Lyons, T. J. Stochastic area for Brownian motion
on the Sierpinski gasket. Ann. Probab. 26, no. 1, 132--148. 1998.

\bibitem{HS} Hebisch, W.; Sikora, A.: A smooth subadditive homogeneous norm
on a homogeneous group. Studia Math. 96, 231-236. 1990.

\bibitem{KS} Karatzas, I.; Shreve, E.:\ Brownian Motion and Stochastic
Calculus. Second Edition. 1991.

\bibitem{L94} Lyons, T.: Differential Equations driven by rough signals.
Math. Res. Letters 1, 451-464. 1994.

\bibitem{Lej} Lejay, A.: Introduction to Rough Paths, S\'{e}minaire de
probabilit\'{e}s, Lecture Notes in Mathematics, Springer-Verlag, volume
XXXVII. 2003.

\bibitem{Le} Levy, P., Processus Stochastiques et Mouvement Brownien.
Gauthier-Villars, 1948.

\bibitem{Ly} Lyons, T.: Differential equations driven by rough signals. Rev.
Mat. Iberoamericana 14, no. 2, 215--310, 1998.

\bibitem{LQ} Lyons, T.; Qian, Z.: System Control and Rough Paths, Oxford
University Press 2002.

\bibitem{LQZ} Ledoux, M.; Qian, Z.; Zhang, T. Large deviations and support
theorem for diffusion processes via rough paths. Stochastic Process. Appl.
102, no. 2, 265--283. 2002

\bibitem{Ma} Malliavin, P: Stochastic analysis. Grundlehren der
Mathematischen Wissenschaften, 313. Springer-Verlag, Berlin, 1997.

\bibitem{Me} Mellouk, M. : Support des diffusions dans les espaces de
Besov-Orlicz. C. R. Acad. Sci. Paris S\'{e}r. I Math. 319, no. 3, 261--266.
1994.

\bibitem{Millet} Millet, A.; Sanz-Sol\'{e}, M.: A simple proof of the
support theorem for diffusion processes. S\'{e}minaire de Probabilit\'{e}s,
XXVIII, 36--48, Lecture Notes in Math., 1583, Springer, Berlin, 1994.

\bibitem{Mo} Montgomery, R.: A tour of subriemannian geometries, their
geodesics and applications. Mathematical Surveys and Monographs, 91.
American Mathematical Society, Providence, RI, 2002.

\bibitem{Ne} Neuenschwander, D.: Probability on the Heisenberg Group,
Springer-LNM. 1996.

\bibitem{Re} Reutenauer C.: Free Lie algebras. London Mathematical Society
Monographs. New Series, 7. Oxford Science Publications, 1993.

\bibitem{RY} Revuz D., Yor N.: Continuous Martingales and Brownian Motion,
Springer. 2001.

\bibitem{Sip} Sipilainen, E.M.: A pathwise view of solutions of stochastic
differential equations. Ph.D. Thesis, University of Edinburgh, 1993.

\bibitem{Ste} Stein, E.M.: Singular integrals and differentiability
properties of functions. Princeton Mathematical Series, No. 30, 1970.

\bibitem{ST} Stroock, D.; Taniguchi, S.: Diffusions as integral curves, or
Stratonovich without It\^{o}. The Dynkin Festschrift, 333--369, Progr.
Probab., 34, 1994.

\bibitem{SV} Stroock, D.; Varadhan, S.R. : Multidimensional diffusion
processes. Springer. 1979.

\bibitem{SV2} Stroock, D.; Varadhan, S. R.: On the support of diffusion
processes with applications to the strong maximum principle. Proceedings of
the Sixth Berkeley Symposium on Mathematical Statistics and Probability,
Vol. III: Probability theory, pp. 333--359, 1972.

\bibitem{VaLN} Varadhan, S.R.S.: Brownian Motion. Lecture notes.

\bibitem{Va} Varadhan, S.R.S.: Large Deviations, SIAM. 1984.

\bibitem{Var} Varadarajan, V.S.: Lie Groups, Lie Algebras and their
presentations, Springer. 1984.

\bibitem{Vi} Victoir, N.: An extension theorem to rough path, (preprint
2003).

\bibitem{Wa} Warner, F.: Foundations of Differentiable Manifolds and Lie
Groups, Springer. 1983.

\bibitem{Yo} Young, L. C.: An inequality of H\"{o}lder type, connected with
Stieltjes integration. Acta Math. 67, 251-282. 1936.
\end{thebibliography}
\end{document}